\documentclass[a4paper,11pt]{article}
\usepackage{latexsym,amsfonts,amsmath}
\usepackage{amssymb}
\usepackage{color}
\usepackage[mathscr]{eucal}
\usepackage{enumerate}
\usepackage{enumitem}

\def\1{\mathbf{1}}
\def\0{\mathbf{0}}
\def\inter{\mathop{\cap}}
\def\NN{\mathbb{N}}
\def\ZZ{\mathbb{Z}}
\def\RR{\mathbb{R}}

\def\XX{\mathbf{X}}
\def\AA{\mathbf{A}}
\def\KK{\mathbf{K}}

\newcommand{\widebar}[1]{\overline{#1}}

\newcommand{\bcal}[1]{\boldsymbol{\mathcal{#1}}}
\newcommand{\bscr}[1]{\boldsymbol{\mathscr{#1}}}
\newcommand{\bfrak}[1]{\boldsymbol{\mathfrak{#1}}}

\def\union{\mathop{\cup}}
\def\inter{\mathop{\cap}}

\def\ds{\displaystyle}
\def\nl{\mbox{} \newline }

\newcounter{hypot}

    \newenvironment{hypot}{\begin{list}
      {\hspace{\labelsep}\bfseries Assumption \Alph{hypot}.}
      {\leftmargin=0pt
       \labelwidth=0cm
       \refstepcounter{hypot}
       \def\makelabel##1{##1}}}{\end{list}}

\newcounter{assump}

\newenvironment{assump}{\begin{list}
      {\hspace{\labelsep} (\Alph{hypot}.\arabic{assump})}
      {\leftmargin=0pt
       \labelwidth=0cm
       \usecounter{assump}
       }}{\end{list}}

\newcounter{cond}
\newenvironment{condit}[1]{\begin{list}
      {\hspace{\labelsep} (#1\arabic{cond})}
      {\leftmargin=0pt
       \labelwidth=0cm
       \usecounter{cond}
       }}{\end{list}}

\begin{document}
\newtheorem{theorem}{Theorem}[section]
\newtheorem{proposition}[theorem]{Proposition}
\newtheorem{lemma}[theorem]{Lemma}
\newtheorem{corollary}[theorem]{Corollary}
\newtheorem{definition}[theorem]{Definition}
\newtheorem{remark}[theorem]{Remark}
\newtheorem{conjecture}[theorem]{Conjecture}
\newtheorem{assumption}[theorem]{Assumption}

\bibliographystyle{plain}

\title{A convex programming approach for discrete-time Markov decision processes under the expected total reward criterion}

\author{ \mbox{ }
F. Dufour \\
\small Institut Polytechnique de Bordeaux \\
\small INRIA Bordeaux Sud Ouest, Team: CQFD\\
\small IMB, Institut de Math\'ematiques de Bordeaux, Universit\'e de Bordeaux, France\\
\small e-mail: francois.dufour@math.u-bordeaux.fr
\and
A. Genadot \\
\small IMB, Institut de Math\'ematiques de Bordeaux, Universit\'e de Bordeaux, France\\
\small INRIA Bordeaux Sud Ouest, Team: CQFD\\
\small e-mail: alexandre.genadot@math.u-bordeaux.fr
%\and
%A. Piunovskiy \\
%\small Department of Mathematical Sciences \\
%\small University of Liverpool  \\
%\small L69 7ZL, Liverpool, United Kingdom \\
%\small e-mail: piunov@liv.ac.uk
}

%\date{Version - final}
\date{}

\maketitle

\begin{abstract} 
In this work, we study discrete-time Markov decision processes (MDPs) under constraints with Borel state and action spaces and where all the performance functions have the same form of the expected total reward (ETR) criterion over the infinite time horizon.
%The main topic of this paper is to study this type of MDPs through a linear programming formulation.
One of our objective is to propose a convex programming formulation for this type of MDPs.
%\textcolor{red}{In particular, it will be shown that the original constrained control problem is equivalent to an associated linear program in the sense that 
%the corresponding values are equal and that an optimal admissible solution for the latter yields an optimal admissible policy for the former and \textit{vice versa}.}
It will be shown that the values of the constrained control problem and the associated convex program coincide and that if there exists an optimal solution to the convex program then there exists a stationary randomized policy which is optimal for the MDP.
%It is important to emphasize that such results are very scarce in the literature.
It will be also shown that in the framework of constrained control problems, the supremum of the expected total rewards over the set of randomized policies is equal to the supremum of the expected total rewards over the set of stationary randomized policies.
We consider standard hypotheses such as the so-called continuity-compactness conditions and a Slater-type condition. Our assumptions are quite weak to deal with cases that have not yet been addressed in the literature.
An example is presented to illustrate our results with respect to those of the literature.
\end{abstract}

{\small
\par\noindent\textbf{Keywords:} Markov decision process, expected total reward criterion, occupation measure, constraints, convex program.
\par\noindent\textbf{AMS 2010 Subject Classification:} 90C40, 60J10, 90C90.}

%%%%%%%%%%%%%%%%%%%%%%%%%%%%%%%%%%%%%%%%%%%%%%%%%%%%%%
\section{Introduction}
\label{sec-1}
We consider a discrete-time Markov decision process with constraints when all the objectives have the same form of the expected total reward over the infinite time horizon.
Markov decision processes are a general family of controlled stochastic processes, which are suitable for the modeling of sequential decision-making problems under uncertainty.
They arise in many applications, such as engineering, medicine, biology, operations research, management science, economics, among others.

Markov decision processes (MDPs) under the expected total reward (ETR) criterion have been extensively studied by using mainly different approaches, see \textit{e.g.} \cite{feinberg02} for a complete and exhaustive survey on that subject and also \cite[Chapter 2]{piunovskiy13} for an analysis of that topic through examples.

When dealing with constraints, the linear/convex programming approach (also called the convex analytic method, see, e.g. \cite{borkar88,borkar02}) has proved to be a very powerful technique for solving MDPs.
It has been extensively studied in the literature and we refer the interested reader to the following works \cite{altman99,borkar88,borkar02,feinberg12,piunovskiy97} and the references therein to get an overview of this technique. The convex programming approach can be applied to a large class of control problems including for example, the finite-horizon and the infinite-horizon discounted-reward problems; see, e.g., \cite{borkar02} for further examples of performance functions. For such criteria, the key idea is to reformulate the original dynamic control problem as an infinite dimensional static optimization problem over a space of finite measures given by the occupation measures of the controlled process. However, it must be emphasized that the expected total reward criterion is an exception where the convex programming formulation may not be suitable except for very specific models. As mentioned in \cite[p. 357-358]{borkar02} and \cite[p. 92-93]{hernandez99}, the ETR criterion is
  very demanding from a technical point of view and yields some important technical difficulties which are basically of two types:
\begin{enumerate}[label={\alph*)}]
\item The first issue is directly related to the question of how to properly formulate a convex program associated with an MDP under the ETR criterion.
Indeed, as described in \cite{borkar02}, the classical and natural approach to formulate a convex program associated to a MDP is to consider as underlying vector space the set of signed finite measures and as variables the occupation measures of the process. However, in the context of the ETR criterion, this approach fails since the occupation measures are not necessarily finite and may take the value infinity. Therefore, the space of finite signed measures is not the appropriate vector space to define the convex program.
%so, the vector space given by the set of finite measure cannot be used to describe the elements of the linear program.
\item An important issue is related to the so-called characteristic equation satisfied by the occupation measures of the process which is of the form:
$$\mu_X (\cdot)= \nu(\cdot) + \int_{X\times A} Q(\cdot |x,a) \mu(dx,da)$$
where $X$ and $A$ are respectively the state and action spaces; $Q$ is the transition probability function of the MDP and $\mu_X$ is the marginal of the measure $\mu$ on $X$.
Indeed, a solution $\mu$ to this equation may not correspond to any occupation measures of the controlled process.
This difficulty makes the analysis of the ETR criterion very involved by using the convex programming approach.
\end{enumerate}

The objective of the current paper is to propose a suitable convex program for MDPs under the ETR criterion.
Our purpose is also to show that the value of the constrained control problem corresponds to the value of an associated convex program and that if there exists an optimal solution to the associated convex program then there exists a stationary randomized policy which is optimal for the MDP.
We consider standard assumptions, the so-called continuity-compactness conditions introduced by Sch\"{a}l in
\cite{schal75,schal83}. These assumptions are of two types, namely conditions (S) and (W). Roughly speaking condition (S) requires the transition kernel to be \textit{strongly continuous}
whereas condition (W) refers to the case where the transition kernel is \textit{weakly continuous}, see, e.g., \cite[p. 367-368]{schal83} for a precise statement of these assumptions. 
We also suppose the existence of a policy in the \textit{interior} of the set of admissible policies. This is the so-called Slater condition.
Conditions (W) and (S) do not play the same role in the sense that when working with condition (W) instead of condition (S) we have to consider an additional hypothesis requiring the transition kernel of the model to be absolutely continuous with respect to a Markov kernel uniformly in the action variables.
Our approach differs from that classically considered in the literature in the sense that the variables of the convex program are not given by the occupation measures of the controlled process but defined on the positive cone of the vector space given by the pair of finite signed stochastic kernels on the action space given the state space.

When compared to the literature, our results appear complementary and our assumptions are rather weak.
The references dealing with the ETR criterion by using the convex programming formulation are very scarce in the literature.
As for our work, the results in \cite{dufour12,dufour13} are concerned with general Borel state and action spaces.
However, it is important to observe that the approach proposed in \cite{dufour12,dufour13} does not correspond to a linear/convex programming formulation of an  MDP under the ETR criterion.
Indeed, the underlying variables of the optimization problem under consideration are given by measures that may take the value infinity and therefore, {this set does not enjoy the structure of a standard vector space.} 
This technical issue aside, the results of the current paper differ significantly from those obtained in \cite{dufour12,dufour13}.
The approach developed in \cite{dufour12} deals with models satisfying condition (W) and strongly relies on the positiveness of the cost functions.
It must be emphasized that the general framework of signed cost functions cannot be addressed with the technique presented in \cite{dufour12}.
In \cite{dufour13}, the model under consideration satisfies condition (S) and it was assumed that the transition kernel is absolutely continuous with respect to a reference probability measure uniformly in the state and action variables. In the present work, we show that this assumption is not needed under condition (S).
It must be also observed that the approach developed in \cite{dufour13} for signed  cost function cannot be applied under condition (W).
In \cite[Chapter 8]{altman99}, the model is transient or absorbing and is restricted to discrete state and action spaces.
Here, we do not impose the MDP to be transient or absorbing.
Another advantage of our approach is to propose a convex programming formulation for constrained MDPs under the ETR criterion with signed reward functions and satisfying condition (W). In this context, such formulation has not been so far investigated in the literature.
It should be also mentioned that in our work we imposed the so-called Slater condition which is not required in \cite{altman99,dufour12,dufour13}.
However, this condition is rather weak and it is a standard assumption in convex optimization problems with constraints, see e.g. \cite{barbu12}.
%\textcolor{blue}{Our approach has the advantage of dealing with cases that have not been yet addressed in the literature such as for example we propose a linear program for constrained MDPs under the ETR criterion with signed reward functions and satisfying condition (W).}

The rest of the paper is organized as follows.
In Section \ref{def-ass}, we present the control problem that will be considered throughout this work.
The assumptions and the convex programming formulation of a constrained discrete-time MDP under the ETR criterion is introduced in Section \ref{sec-Assumptions-LP-formulation}.
Important properties of the convex program as well as the constrained control problem are established in Section \ref{sec-Preliminary-Results}.
Our main results are presented in Section \ref{LP-formulation=Control-Problem} showing that the original control problem is equivalent to the convex program.
Section \ref{sec-example} is dedicated to the presentation of an example illustrating our results.
Finally, a technical result used in Section \ref{sec-Preliminary-Results} is derived in an appendix.

%%%%%%%%%%%%%%%%%%%%%%%%%%%%%%%%%%%%%%%%%%%%%%%%%%%%%%

%%%%%%%%%%%%%%%%%%%%%%%%%%%%%%%%%%%%%%%%%%%%%%%%%%%%%%
\section{Description of the control problem}
\label{def-ass}
%%%%%%%%%%%%%%%%%%%%%%%%%%%%%%%%%%%%%%%%%%%%%%%%%%%%%%
\setcounter{equation}{0}
The main goal of this section is to introduce the notation, the parameters defining the model, and to present the construction of the controlled process.

\subsection{Notation and terminology}\label{sec-notation}
The following basic notation will be used in the forthcoming.

The set of integers is denoted by $\ZZ$ and $\NN$ corresponds to the non-negative integers, that is, $\NN=\{0,1,2,\ldots\}$. The set of real numbers is given by $\RR$.
For any subset $\mathbb{D}$ of $\RR$, $\mathbb{D}^{*}$ denotes $\mathbb{D}\setminus \{0\}$ and $\mathbb{D}_{+}=\{d\in \mathbb{D}: d\geq 0\}$.
We write $\NN_{p}$ for $\{1,\ldots,p\}$ with $p\in \NN^{*}$, $\widebar{\RR}$ is the set of extended real numbers, that is, $\RR\cup\{-\infty,+\infty\}$ and $\widebar{\RR}_{+}=\RR_{+}\cup\{+\infty\}$.
Given $x$ and $y$ in the Euclidean space $\RR^{n}$,
let $\langle x,y\rangle$ be the usual inner product of $x$ and $y$. By $|x|=\langle x,x\rangle^{1/2}$ we will denote the norm of $x\in\RR^n$. Let $\0_{n}$ be the element of $\RR^n$ with all components equal to zero.
If $\theta_1$ and $\theta_2$ are in $\RR^{n}$, we shall write $\theta_1\ge\theta_2$ when all the components of $\theta_1$ are greater than or equal to the corresponding components of~$\theta_2$.

Let $X$ be a metric space and denote by $\bfrak{B}(X)$ its associated Borel $\sigma$-algebra.
We use the symbol $f^{+}$ (respectively $f^{-}$) to denote the positive part (respectively, negative part) of a function
$f:X\rightarrow\widebar{\RR}$.
The function $\mathcal{I}_{\infty}$ is the function whose values are constant and equal to $+\infty$.
If $X$ is a metric space, $\bscr{M}(X)$ denotes the set of real-valued measurable functions defined on $X$.
Furthermore, $\bcal{C}(X)$ is the space of real-valued bounded continuous functions defined on $X$.
The term measure will always refer to a countably additive, $\widebar{\RR}_{+}$-valued set function.
The set of measures defined on $(X,\bfrak{B}(X))$ is denoted by $\bcal{M}(X)$ and the set of probability measures on $(X,\bfrak{B}(X))$ by $\boldsymbol{\mathcal{P}}(X)$.
For $\mu\in \bcal{M}(X)$ and a positive function $h$ in $\bscr{M}(X)$, $\mu(h)=\int_{X} h(x) \mu(dx)$ and for $g\in \bscr{M}(X)$, $\mu(g)$ is defined by $\mu(g^+)-\mu(g^-)$ where by convention $(+\infty)-(+\infty)=-\infty$.
Consider two metric spaces $X$ and $Y$.
If $\mu$ is a measure on $X\times Y$ then $\mu_{X}$ denotes the marginal of the measure $\mu$ on $X$. 
A kernel $K$ on $X$ given $Y$ is a $\widebar{\RR}_{+}$-valued mapping defined on $\bfrak{B}(X)\times Y$ such that
for any $y\in Y$, $K(\cdot |y)\in \bcal{M}(X)$ and for any $\Lambda\in\bfrak{B}(X)$, $K(\Lambda|\cdot)$ is a measurable function defined on $Y$.
A kernel $K$ on $X$ given $Y$ is said to be finite if $K(X |y)\in \RR_{+}$ for any $y\in Y$.
The set of finite kernels on $X$ given $Y$ is denoted $\bcal{K}(X|Y)$. A stochastic (or Markov) kernel $K$ on $X$ given $Y$ is a kernel in $\bcal{K}(X|Y)$ satisfying
$K(X |y)=1$ for any $y\in Y$. The set of stochastic kernels on $X$ given $Y$ will be denoted by $\bcal{P}(X|Y)$.
Let $Q$ be a stochastic kernel on $X$ given $Y$, then, for a function $v:X\rightarrow\widebar{\RR}$, we define $Qv:Y\rightarrow\widebar{\RR}$ as
$$Qv(y):=\int_{X}v^{+}(x)Q(dx|y)-\int_{X}v^{-}(x)Q(dx|y),$$
provided that $v$ is quasi-integrable with respect to the probability measure $Q(\cdot |y)$ for any $y\in Y$.
For a measure $\mu$ on $Y$, we denote by $\mu Q$ the measure $\ds \int_{Y} Q(\cdot|y)\mu(dy)$ on $X$.
%For a finite kernel $K$ on $X$ given $Y$ and a probability measure $q$ on $Y$, the product measure $K\cdot q$ on $X\times Y$ is defined by $\ds K\cdot q (\Gamma) = \int_{\Gamma} K(dx|y) q(dy)$.
%Furthermore, $\bcal{C}(X)$ (respectively, $\bfrak{U}(X)$) is the space of real-valued continuous (respectively, upper-semicontinuous) function defined $X$.
%If $\gamma\in \bcal{M}(X)$ and $v\in \bscr{M}_{\RR^{p}}(X)$, then the (componentwise) integral of $v$ with respect to $\gamma$, provided that it is well defined, will be denoted by $\ds \gamma(v)=\int_{X} v d\gamma\in\widebar{\RR}^{p}$.

\subsection{The control model.}
Let us consider the stationary model 
\begin{equation}
\big(\mathbf{X},\mathbf{A},\{\mathbf{A}(x): x\in \mathbf{X}\},Q,r,c,\theta_{*},\nu\big) \label{4tuple}
\end{equation}
consisting of:
\begin{enumerate}
\item [(a)] A Borel space $\mathbf{X}$ (that is, a Borel subset of a complete and separable metric space), which is the state space.
\item [(b)] A Borel space $\mathbf{A}$, representing the control or action set.
\item [(c)] A family $\{\mathbf{A}(x): x\in \mathbf{X}\}$ of non-empty measurable subsets of $\mathbf{A}$, where $\mathbf{A}(x)$ is the set of feasible controls or actions when the system is in state $x\in \mathbf{X}$. We suppose that
\begin{eqnarray*}
\mathbf{K} :=\{(x,a)\in \mathbf{X}\times \mathbf{A}: a\in \mathbf{A}(x)\}
\end{eqnarray*}
is a measurable subset of $\mathbf{X}\times \mathbf{A}$.
There exists a measurable map $\vartheta : \mathbf{X} \rightarrow \mathbf{A}$ with $\vartheta(x)\in\mathbf{A}(x)$.
For notational convenience, we introduce recursively the set $\mathbf{H}_{t}$ of histories up to time $t\in \NN^{*}$ by defining $\mathbf{H}_{1}=\mathbf{X}$ and
$\mathbf{H}_{t}=\mathbf{K}^{t-1}\times \mathbf{X}$ for $t\geq 2$.
\item [(d)] A stochastic kernel $Q$ on $\mathbf{X}$ given $\mathbf{K}$, which stands for the transition probability function. 
\item [(e)] The one-step reward function is given by a measurable function $r:\mathbf{K}\rightarrow \RR$.
\item [(f)] For $i\in\NN_{q}$, the measurable mappings $c_{i}:\mathbf{K}\rightarrow \RR$ are the one-step constraint functions.
\item [(g)] The constraint limits are real numbers given by $\theta^{*}=\big\{\theta^{*}_{i}\big\}_{i\in\NN_{q}}$.
\item [(h)] Finally, the initial distribution is $\nu\in\bcal{P}(\XX)$.
\end{enumerate}
A control policy (a policy, for short) is a sequence $\pi=\{\pi_{t}\}_{t\in\NN^{*}}$ of stochastic kernels $\pi_{t}$ on~$\mathbf{A}$ given $\mathbf{H}_{t}$ such that $\pi_{t}(\mathbf{A}(x_{t})|h_{t})=1$ for any $h_{t}=(x_{1},a_{1},\ldots,x_{t})\in \mathbf{H}_{t}$.
Let $\Pi$ be the set of all policies.
A policy $\pi=\{\pi_{t}\}_{t\in\NN^{*}}\in \Pi$ is called a stationary randomized policy if there exists a stochastic kernel $\boldsymbol{\varphi}$ on $\AA$ given $\XX$ satisfying $\boldsymbol{\varphi}(\AA(x)|x)=1$
for any $x\in \XX$ and $\pi_{t}(\cdot |h_{t})=\boldsymbol{\varphi}(\cdot |x_{t})$ for any $h_{t}=(x_{1},a_{1},\ldots,x_{t})\in \mathbf{H}_{t}$ and $t\in \NN^{*}$.
{In such as case, we will write $\boldsymbol{\varphi}$ instead of $\pi$ to emphasize that the corresponding stationary randomized policy $\pi$ is generated by $\boldsymbol{\varphi}$.}
Let $\Pi_{s}$ be the set of all stationary randomized policies.

\bigskip

To state the optimal control problem we are concerned with, we introduce the canonical space $(\Omega,\mathcal{F})$ consisting of the set of sample paths $\Omega=(\XX\times\AA)^{\infty}$
and the associated product $\sigma$-algebra $\mathcal{F}$. The projection from $\Omega$ to the state space and the action space at time $t$ are denoted by $X_{t}$ and $A_{t}$. That is, for
\[\omega=(y_{1},b_{1},\ldots,y_{t},b_{t}\ldots)\in \Omega \quad\hbox{we have}\quad X_{t}(\omega)=y_{t}\;\;\hbox{and}\;\; A_{t}(\omega)=b_{t}\] for $t\in \NN^{*}$.
Consequently, $\{X_{t}\}_{t\in \NN^{*}}$ is the state process and $\{A_{t}\}_{t\in \NN^{*}}$ is the control process.
It is a well known result that for every policy $\pi \in \Pi$ and any initial probability measure $\nu$ on $(\mathbf{X},\bfrak{B}(\mathbf{X}))$ there exists a unique probability
measure $\mathbb{P}_{\nu}^{\pi}$ on $(\Omega,\mathcal{F})$
such that $ \mathbb{P}_{\nu}^{\pi} (\KK^{\infty})=1$ and
$$\mathbb{P}_{\nu}^{\pi}(X_{1}\in B)=\nu(B), \quad \text{ for } B\in \bfrak{B}(\mathbf{X}),$$
$$\mathbb{P}_{\nu}^{\pi}(X_{t+1}\in C|\sigma\{X_{1},\ldots,X_{t},A_{t}\})=Q(C|X_{t},A_{t}) \quad \text{ for } C\in \bfrak{B}(\mathbf{X}),$$ 
$$ \mathbb{P}_{\nu}^{\pi}(A_{t}\in D| \sigma\{X_{1},\ldots,X_{t-1},A_{t-1},X_{t}\})=\pi_{t}(D|X_{1},\ldots,X_{t-1},A_{t-1},X_{t}) \quad \text{ for } D\in \bfrak{B}(\mathbf{A}),$$
$\mathbb{P}_{\nu}^{\pi}-a.s.$, for any $t\in\NN^{*}$.

The expectation with respect to  $\mathbb{P}_{\nu}^{\pi}$ is denoted by $\mathbb{E}_{\nu}^{\pi}$.
The so-called \textit{occupation measure} generated by a policy $\pi\in\Pi$, denoted by $\mu^{\pi}$, is defined by
$$\mu^{\pi}(\Gamma) = \sum_{t=1}^{\infty} \mathbb{P}_{\nu}^{\pi} \big( (X_{t},A_{t})\in \Gamma \big)$$
for any $\Gamma\in \bfrak{B}(\mathbf{X}\times \mathbf{A})$.
Denote by $\bcal{O}$ (respectively, $\bcal{O}_{s}$) the set of occupation measures generated by randomized (respectively, stationary) policies.

\paragraph{Statement of the control problem.}
\nl
For $h\in\bscr{M}(\KK)$ and $\pi\in\Pi$,  define $\mathcal{J}_{\nu}(h,\pi)$ by
\begin{eqnarray*}
\mathcal{J}_{\nu}(h,\pi)=\sum_{t=1}^{\infty} \mathbb{E}_{\nu}^{\pi} \big[h^{+}(X_{t},A_{t})\big] - \sum_{t=1}^{\infty} \mathbb{E}_{\nu}^{\pi} \big[h^{-}(X_{t},A_{t})\big] 
\end{eqnarray*}
where by convention $(+\infty)-(+\infty)=-\infty$. In fact, assumptions will be introduced in the next section to avoid dealing with such cases. Observe that $\mathcal{J}_{\nu}(h,\pi)$ can be written equivalently in terms of the occupation measure generated by the policy $\pi\in\Pi$ as follows
\begin{eqnarray*}
\mathcal{J}_{\nu}(h,\pi)=\mu^{\pi}(h).
\end{eqnarray*}
In this paper, we will repeatedly use this equality without mentioning it.

\begin{definition}
\label{Constrained-control}
A policy $\pi\in \Pi$ is said to be admissible if $\mathcal{J}_{\nu}(c_{i},\pi)\geq \theta^{*}_{i}$ for $i\in \NN_{q}$. The set of admissible policies will be denoted by
$\Pi_{\theta^{*}}$.
The optimal control problem we consider consists in maximizing the expected reward $\mathcal{J}_{\nu}(r,\pi)$ over the set of admissible policies $\pi\in\Pi_{\theta^{*}}$.
The value associated to this constrained control problem is given by $\sup\big\{ \mathcal{J}_{\nu}(r,\pi) : \pi\in\Pi_{\theta^{*}} \big\} $.
A policy $\hat{\pi}\in \Pi$ is optimal if $\hat{\pi}\in \Pi_{\theta^{*}}$ and $\mathcal{J}_{\nu}(r,\hat{\pi}) =\sup\big\{ \mathcal{J}_{\nu}(r,\pi) : \pi\in\Pi_{\theta^{*}} \big\} $.
\end{definition}

%%%%%%%%%%%%%%%%%%%%%%%%%%%%%%%
\section{Assumptions and the convex programming formulation}
\label{sec-Assumptions-LP-formulation}
%%%%%%%%%%%%%%%%%%%%%%%%%%%%%%%
The objective of this section is both to list the assumptions we will use in this work and to introduce the convex program associated with the control problem we presented in the previous section.
In this work, we deal with MDPs satisfying the so-called Conditions (W) or (S) which are standard hypotheses of the literature, see for example  \cite{schal75}. 

\bigskip

\noindent \textbf{Condition (W)}
\begin{condit}{W}
\item \label{Hyp-Multi-function-upper-semicontinuous} For any $x\in \XX$, the action set $\AA(x)$ is compact and the multifunction from $\XX$ to~$\AA$ defined by
$x\rightarrow \AA(x)$ is upper-semicontinuous.
\item \label{Hyp-Q-weakly-continuous} For any $f\in \bcal{C}(\XX)$, $Qf$ is continuous on $\KK$.
\item \label{Hyp-W-reward-constraint-upper-semicontinuous} The reward $r$ and the constraint $c_{i}$ for $i\in\NN_{q}$
are upper-semicontinuous on $\KK$.
\end{condit}

\noindent\textbf{Condition (S)}
\begin{condit}{S}
\item \label{Hyp-Compactness} For any $x\in \XX$, $\AA(x)$ is compact.
\item \label{Hyp-Q-strongly-continuous} For any $x\in \XX$ and $\Lambda\in \bfrak{B}(\XX)$, $Q(\Lambda |x,\cdot)$ is continuous on $\AA(x)$.
\item \label{Hyp-S-reward-constraint-upper-semicontinuous} For any $x\in \XX$, the reward $r(x,\cdot)$ and the constraint $c_{i}(x,\cdot)$ for $i\in\NN_{q}$
are upper-semicontinuous on $\AA(x)$.
\end{condit}

In order to introduce the convex program associated to an MDP under the ETR criterion, we need to make some hypotheses.
First, it is assumed that the transition kernel $Q$ of the MDP under consideration is absolutely continuous with respect to a Markov kernel $P$ (see Assumption
\ref{Hyp-Transition-kernel-absolute-continuity}).
This hypothesis is rather weak and is satisfied in a large number of practical cases as discussed in the remark below. 

\begin{hypot}
\label{Hyp-Transition-kernel-absolute-continuity}
\item There exists $P\in \bcal{P}(\XX| \XX)$ satisfying $Q(\cdot |x,a)\ll P(\cdot |x)$ for any $(x,a)\in\mathbf{K}$. 
Associated to the kernel $P$, $p$  will denote the probability measure on $\XX$ defined by
\begin{eqnarray}
p(dx) = \sum_{k\in \NN} \frac{1}{2^{k+1}} \nu P^{k}(dx).
\label{Def-p}
\end{eqnarray}
\end{hypot}

\begin{remark}
\label{Discussion-Hyp-Transition-kernel-absolute-continuity}
\begin{enumerate}
\item In Lemma \ref{Existence-P-strongly-continuous} below, it is shown that under Conditions (S\ref{Hyp-Compactness}) and (S\ref{Hyp-Q-strongly-continuous}), Assumption \ref{Hyp-Transition-kernel-absolute-continuity} is satisfied.
\item If the sets of feasible actions are countable, that is $\mathbf{A}(x)=\{a_k(x)\}_{k\in\NN^{*}}$ where for any $k\in \NN^{*}$ $a_{k}$ is a measurable function from $\mathbf{X}$ to
$\mathbf{A}$ then Assumption \ref{Hyp-Transition-kernel-absolute-continuity} is satisfied for $P$ defined by
$$ P(dy|x)=\sum_{k\in \NN^{*}} \frac{1}{2^{k}} Q(dy|x,a_{k}(x)), $$
for any $x\in\mathbf{X}$. 
\item If $Q(\cdot |x,a)\ll q(\cdot)$ for any $(x,a)\in\mathbf{K}$ then clearly Assumption \ref{Hyp-Transition-kernel-absolute-continuity} is satisfied. This condition corresponds to the main hypothesis used in \cite{dufour13}. It is of course less general than Assumption \ref{Hyp-Transition-kernel-absolute-continuity} but it is naturally satisfied for a large class of practical systems.
Indeed,  in many applications, the evolution of an MDP is specified by a discrete-time equation of the form
$x_{t+1}=F(x_{t},a_{t})+\xi_{t}$ where $F$ is an $\RR^{n}$-valued measurable mapping defined on $\RR^{n}\times A$ and $(\xi_{t})_{t\in \NN^{*}}$ is an independent and identically distributed sequence of random variables with density $\alpha$ with respect to the Lebesgue measure on $\bfrak{B}(\RR^{n})$. By using the change of variable formula,
we obtain that $\ds Q(A|x,a)=\int_{A}\alpha(y-F(x,a)) dy$ showing that 
$Q(\cdot |x,a)\ll q(\cdot)$ for any $(x,a)\in\mathbf{K}$ is satisfied for $q$ defined for example by the standard normal distribution on $\bfrak{B}(\RR^{n})$.
\nl
Observe also that when $\mathbf{X}$ is finite or countable, $Q(\cdot |x,a)\ll q(\cdot)$ for any $(x,a)\in\mathbf{K}$ is satisfied when $q$ is given for example by a geometric distribution.
\end{enumerate}
\end{remark}

\begin{lemma}
\label{Existence-P-strongly-continuous}
Conditions (S\ref{Hyp-Compactness}) and (S\ref{Hyp-Q-strongly-continuous}) imply Assumption \ref{Hyp-Transition-kernel-absolute-continuity}, that is, $Q\ll P$ with
$P\in \bcal{P}(\XX|\XX)$  given by
\begin{eqnarray}
P(dy|x)=\sum_{k\in \NN^{*}}  \frac{1}{2^{k}} Q(dy|x,\xi_{k}(x))
\label{Def-P-strongly-continuous}
\end{eqnarray}
where $\{\xi_{k}\}_{k\in\NN^{*}}$ is a sequence of measurable selectors from the multifunction defined from $\XX$ to~$\AA$ by $x\rightarrow \AA(x)$ and
satisfying $\AA(x)= \widebar{\{\xi_{n}(x): n\in\NN^{*}\}}$ for any $x\in\XX$.
\end{lemma}
\textbf{Proof:} The multifunction $\bfrak{A}$ from $\XX$ to~$\AA$ defined by $x\rightarrow \AA(x)$ is by assumption Borel measurable and so, weakly measurable.
From (S\ref{Hyp-Compactness}), Corollary 18.15 in \cite{aliprantis06} gives the existence of a sequence $\{\xi_{n}\}_{n\in\NN^{*}}$ of measurable selectors from the multifunction $\bfrak{A}$
satisfying $\AA(x)= \widebar{\{\xi_{n}(x): n\in\NN^{*}\}}$ for any $x\in\XX$.
Now by using (S\ref{Hyp-Q-strongly-continuous}), we obtain that
$Q(dy|x,a)\ll P(dy|x)$ for any $(x,a)\in\KK$ for the Markov kernel $P$ defined by  \eqref{Def-P-strongly-continuous}.
\hfill$\Box$
\begin{remark}
The previous proof is an extension of an argument used in the proof of Theorem 1 in \cite[p. 183]{nowak88}. 
\end{remark}

In the next definition, we introduce the set of feasible variables. It will be shown below that it is a convex subset of the vector space of finite signed kernels on $\AA$ given $\XX$.
\begin{definition}
%\label{Variables-Linear-Program}
Suppose Assumption \ref{Hyp-Transition-kernel-absolute-continuity} holds and let $p$ be the measure introduced in \eqref{Def-p}.
\begin{itemize}
\item For $\Phi=(\varphi^{\infty},\varphi^{*})\in\bcal{K}(\AA| \XX)^{2}$, $\eta^{\Phi}$ will denote the measure in $\bcal{M}(\XX\times\AA)$ given by
\begin{eqnarray}
\eta^{\Phi}(dx,da)= \mathcal{I}_{\infty}(x) \varphi^{\infty}(da|x) p(dx) + \varphi^{*}(da|x) p(dx),
\label{Def-eta-Phi}
\end{eqnarray}
recalling that $\mathcal{I}_{\infty}$ is constant function equal to infinity.
\item Consider $\bcal{K}_{p}$ as the set of
$\Phi=(\varphi^{\infty},\varphi^{*})\in \bcal{K}(\AA| \XX)^{2}$ satisfying
$$\varphi^{\infty}(\AA|x)+\varphi^{*}(\AA|x)>0,$$
$$\varphi^{\infty}(\AA(x)^{c}|x)+\varphi^{*}(\AA(x)^{c}|x)=0,$$
and
$$\eta^{\Phi}_{\XX}= \nu + \eta^{\Phi} Q.$$
Any $\Phi\in \bcal{K}_{p}$ induces a measure $\eta^{\Phi}$ that will be called the \textit{$\bcal{K}_{p}$-measure} generated by $\Phi$.
$\bcal{K}_{p}$ is called the set of feasible variables.
\end{itemize}
\end{definition}

\begin{remark}
Observe first that $\alpha\Phi_{1}+(1-\alpha)\Phi_{2}\in \bcal{K}_{p}$ and in particular, 
\begin{eqnarray}
\label{eta-lineaire}
\eta^{\alpha\Phi_{1}+(1-\alpha)\Phi_{2}}=\alpha \eta^{\Phi_{1}}+ (1-\alpha) \eta^{\Phi_{2}},
\end{eqnarray}
for any $\alpha\in[0,1]$ and $(\Phi_{1},\Phi_{2})\in \bcal{K}_{p}^{2}$. Therefore, $\bcal{K}_{p}$ is a convex subset of the vector space of signed finite kernel on $\AA$ given $\XX$.
\end{remark}

\begin{definition}
\label{Induced-policy}
Let $\Phi=(\varphi^{\infty},\varphi^{*})\in \bcal{K}_{p}$. Introduce the kernel $\varphi_{\Phi}$ on $\AA$ given $\XX$ defined by
\begin{eqnarray}
\label{Def-varphi-Phi}
\varphi_{\Phi}(da |x)= I_{\bcal{E}_{\Phi}^{c}}(x) \frac{\varphi^{\infty}(da |x)}{\varphi^{\infty}(\AA |x)} + I_{\bcal{E}_{\Phi}}(x) \frac{\varphi^{*}(da |x)}{\varphi^{*}(\AA |x)}.
\end{eqnarray}
where
\begin{eqnarray}
\label{Def-E-Phi}
\bcal{E}_{\Phi}=\{x\in \XX : \varphi^{\infty}(\AA |x)=0\}.
\end{eqnarray}
Observe that $\varphi_{\Phi}$ is a stochastic kernel satisfying $\varphi_{\Phi}(\AA(x) |x)=1$ for any $x\in \XX$.
The stationary randomized policy $\varphi_{\Phi}$ will be called the policy induced by $\Phi$.
\end{definition}

We will also need the following technical hypothesis:
\begin{hypot}
\label{Hyp-finiteness}
\item \mbox{ }
\begin{assump}
\item \label{Hyp-finiteness-linear}
$\ds  \sup \big\{ \eta^{\Phi} (r^{+}) : \Phi\in\bcal{K}_{p} \big\}$ and $\ds \sup\big\{ \eta^{\Phi}(c^{+}_{i}) : \Phi\in\bcal{K}_{p} \big\}<+\infty$ for any $i\in\NN_{q}$.
\item \label{Hyp-negative-costs-finiteness}  $\mu(r^{-})<+\infty$ and $\mu(c^{-}_{i})<+\infty$ for any $\mu\in\bcal{O}$, $i\in\NN_{q}$.
\end{assump}
\end{hypot}
This hypothesis is comparable to Assumption (A2) introduced in \cite[p. 847]{dufour13}.
Assumption \ref{Hyp-finiteness-linear} essentially imposes that the values of the unconstrained convex programs associated to a reward function given by either $r$ or $c_{i}$ for $i\in\NN_{q}$ are different from $+\infty$ while Assumption \ref{Hyp-negative-costs-finiteness} ensure that
the performance criteria associated to the reward $r$ and the constraints $c_{i}$ for $i\in\NN_{q}$ are not equal $-\infty$. %These are rather weak assumptions. 
In particular, Assumption \ref{Hyp-finiteness-linear} will be used to introduce the linear program.
\begin{definition}
\label{Def-Linear-Program}
Suppose Assumptions \ref{Hyp-Transition-kernel-absolute-continuity} and {\ref{Hyp-finiteness-linear}} hold.
The convex program, denoted by $\bcal{KP}_{p}$, consists in maximizing $\eta^{\Phi}(r)$ over $\Phi\in\bcal{K}_{p}$ subject to $\eta^{\Phi}(c_{i}) \geq \theta^{*}_{i}$ for any $i\in\NN_{q}$. The \textit{value} of the convex program is given by
\begin{align}
\sup\big\{ \eta^{\Phi}(r) : \Phi\in\bcal{K}_{p} \text{ and } \eta^{\Phi}(c_{i}) \geq \theta^{*}_{i} \text{ for } i\in\NN_{q} \big\}.
\label{Def-Linear-Program-eq1}
\end{align}
A variable $\hat{\Phi} \in\bcal{K}_{p}$ is said to be an optimal solution to the convex program $\bcal{KP}_{p}$ if 
$$\eta^{\hat{\Phi}}(r) = \sup\big\{ \eta^{\Phi}(r) : \Phi\in\bcal{K}_{p} \text{ and } \eta^{\Phi}(c_{i}) \geq \theta^{*}_{i} \text{ for } i\in\NN_{q} \big\}$$
and $\eta^{\hat{\Phi}}(c_{i}) \geq \theta^{*}_{i}$ for any $i\in\NN_{q}$.
\end{definition}

\begin{remark}
Let $h$ be a function given by either $r$ or $c_{i}$ for $i\in\NN_{q}$.
From Assumption \ref{Hyp-finiteness-linear}, it follows that $\alpha \eta^{\Phi_{1}}(h)+ (1-\alpha) \eta^{\Phi_{2}}(h)$ is well defined for any $\alpha\in[0,1]$ and $(\Phi_{1},\Phi_{2})\in \bcal{K}_{p}^{2}$.
Therefore, we obtain from equation \eqref{eta-lineaire} that
$$\eta^{\alpha\Phi_{1}+(1-\alpha)\Phi_{2}}(h)=\alpha \eta^{\Phi_{1}}(h)+ (1-\alpha) \eta^{\Phi_{2}}(h)$$
for any $\alpha\in[0,1]$ and $(\Phi_{1},\Phi_{2})\in \bcal{K}_{p}^{2}$.
{This implies that the mathematical program defined in \eqref{Def-Linear-Program-eq1} is indeed a convex program.}
In \text{\cite[p. 153]{barbu12}}, a convex program is written in terms of an infimum.
The $\bcal{KP}_{p}$ program introduced in Definition \ref{Def-Linear-Program} can be equivalently written in terms of an infimum by changing the sign of the reward function.
We prefer to keep this setting to deal with an MDP under a reward optimization criterion.
\end{remark}

\bigskip
\noindent Finally, we introduce an additional standard hypothesis:
\bigskip

\noindent {\textbf{The Slater condition}
\begin{condit}{Slater}
\item[] There exists $\mu^{*}\in\bcal{O}$ such that $\theta^{*}_{i}<\mu^{*}(c_{i})$ for any $i\in\NN_{q}$.
\end{condit}}

%%%%%%%%%%%%%%%%%%%%%%%%%%%%%%%
\section{Preliminary results}
\label{sec-Preliminary-Results}
%%%%%%%%%%%%%%%%%%%%%%%%%%%%%%%
The main goal of this section is to establish several properties of the constrained control problem as well as properties of the convex program.
%%%%%%%%%%%%%%%
\subsection{Properties of the convex program}
\label{subsec-property-linear-prog}
%%%%%%%%%%%%%%%
In this subsection, we will show in Lemma \ref{strategy-stationary-randomized} that for any stationary randomized policy $\pi\in\Pi_{s}$ there exists $\Phi\in \bcal{K}_{p}$ such that
the $\bcal{K}_{p}$-measure generated by $\Phi$ is equal to the occupation measure generated by the stationary randomized policy $\pi$.
An important result which is a cornerstone of the paper is presented at the end of this subsection. It can be roughly stated as follows: for any feasible variable $\Phi\in\bcal{K}_{p}$ of the convex program, the reward $\mathcal{J}_{\nu}(h,\varphi_{\Phi})$
associated to the stationary randomized policy $\varphi_{\Phi}\in \Pi_{s}$
is greater than $\eta^{\Phi}(h)$ for specific functions $h$ that will be discussed in Theorem \ref{Linear-prog=stationary-optimality}.
To get these results, we first need to establish that the occupation measures of the controlled process have a special structure, that is, the marginal on $\mathbf{X}$ of any occupation measure is absolutely continuous with respect to the probability measure $p$ introduced in Assumption \ref{Hyp-Transition-kernel-absolute-continuity}.
\begin{lemma}
\label{Transition-kernel-absolute-continuity}
Suppose Assumption \ref{Hyp-Transition-kernel-absolute-continuity} holds. Then for any $\mu\in\bcal{O}$, 
\begin{eqnarray}
\mu_{\XX}(dx) \ll p(dx) 
\end{eqnarray}
where $p\in\bcal{P}(\XX)$ is defined in \eqref{Def-p}.
\end{lemma}
\textbf{Proof:}  For any $\mu\in\bcal{O}$, it can be easily shown from Lemma 9.4.3 in \cite{hernandez96} the existence of an increasing sequence of finite measures $\{\mu_{k}\}_{k\in \NN^{*}}$ on $\XX$  and a sequence of stochastic kernels $\{\varphi_{k}\}_{k\in \NN^{*}}$ on $\AA$ given $\XX$ satisfying $\varphi_{k}(\AA(x)|x)=1$ and
\begin{eqnarray}
\lim_{k\rightarrow \infty} \mu_{k}(\Lambda)=\mu_{\XX}(\Lambda)
\label{Transition-eq1}
\end{eqnarray}
and
\begin{eqnarray}
\mu_{k+1}(\Lambda)=\nu(\Lambda)+\int_{\XX} \int_{\AA} Q(\Lambda|x,a)\varphi_{k}(da|x) \mu_{k}(dx)
\label{Transition-eq2}
\end{eqnarray}
for $\Lambda\in \bfrak{B}(\XX)$, $k\in \NN^{*}$ and $\mu_{1}=\nu$.
Let us show by induction that $\mu_{k}\ll p$ for any $k\in \NN^{*}$. We have clearly $\mu_{1}\ll p$. Assume that $\mu_{k}\ll p$.
Observe that $\ds \int_{\AA} Q(\cdot|x,a)\varphi_{k}(da|x) \ll P(\cdot |x)$ for any $x\in \XX$ implying that
$$\int_{\XX} \int_{\AA} Q(\cdot|x,a)\varphi_{k}(da|x) \mu_{k}(dx)\ll \int_{\XX} P(\cdot|x) p(dx)$$
and so, combining \eqref{Def-p} and \eqref{Transition-eq2} we have $\mu_{k+1}\ll p$. We obtain the result by using \eqref{Transition-eq1}.
\hfill$\Box$

\bigskip

As a consequence, we can show that the set of the $\bcal{K}_{p}$-measures contains the occupation mesures generated by the stationary randomized policies.
\begin{lemma}
\label{strategy-stationary-randomized}
Suppose Assumption \ref{Hyp-Transition-kernel-absolute-continuity} holds. 
For any $\pi\in\Pi_{s}$, there exists $\Phi\in  \bcal{K}_{p}$
such that $$\mu^{\pi}=\eta^{\Phi}.$$ 
\end{lemma}
\textbf{Proof:} 
Let $\pi\in \Pi_{s}$.
Clearly, the increasing sequence $\{\mu^{\pi}_{t}\}_{t\in \NN^{*}}$ of finite measures defined on $\XX\times\AA$ by
$$\mu^{\pi}_{t} (\Gamma)=\sum_{k=1}^{t} \mathbb{P}_{\nu}^{\pi} \big( (X_{k},A_{k})\in \Gamma \big)$$
for any $\Gamma\in  \bfrak{B}(\XX\times\AA)$ converges to $\mu^{\pi}$.
From Lemma \ref{Transition-kernel-absolute-continuity}, there exists a sequence of increasing measurable $\RR_{+}$-valued functions
$\{\mathcal{D}_{t}\}_{t\in \NN^{*}}$ defined on $\XX$ such that 
$\ds \sum_{k=1}^{t} \mathbb{P}_{\nu}^{\pi} ( X_{k}\in \Lambda ) = \int_{\Lambda} \mathcal{D}_{t}(x) p(dx)$
for $\Lambda\in \bfrak{B}(\mathbf{X})$ and so,  $\mu^{\pi}_{t}(dx,da)= \mathcal{D}_{t}(x) \pi(da|x)  p(dx)$.
Therefore, 
\begin{eqnarray*}
\mu^{\pi}(dx,da) & = & \mathcal{D}(x) \pi(da|x) p(dx) \nonumber \\
& = & \mathcal{I}_{\infty}(x) I_{\{\mathcal{D}(x)=\infty\}} \pi(da|x) p(dx)
+ \mathcal{D}(x) I_{\{\mathcal{D}(x)<\infty\}} \pi(da|x) p(dx)
\end{eqnarray*}
where $\mathcal{D}(x)=\lim_{t\rightarrow\infty} \mathcal{D}_{t}(x)$.
Consequently, $\Phi=(\varphi^{\infty},\varphi^{*})$ defined by $\varphi^{\infty}(da|x)=I_{\{\mathcal{D}(x)=\infty\}} \pi(da|x)$ and
$\varphi^{*}(da|x)=\mathcal{D}(x) I_{\{\mathcal{D}(x)<\infty\}} \pi(da|x)$ belongs to $\bcal{K}_{p}$
since $\mu^{\pi}_{\XX}= \nu + \mu^{\pi} Q$.
\hfill$\Box$

\bigskip

The following result is in a way a converse of the previous one. It is a key result in our work.
Roughly speaking, it states that for any feasible variable $\Phi\in\bcal{K}_{p}$ of the convex program, the reward $\mathcal{J}_{\nu}(h,\varphi_{\Phi})$ associated to the stationary randomized policy $\varphi_{\Phi}\in \Pi_{s}$ is greater than $\eta^{\Phi}(h)$ for specific functions $h$ described below.
\begin{theorem}
\label{Linear-prog=stationary-optimality}
Suppose that Assumption \ref{Hyp-Transition-kernel-absolute-continuity} holds.
For any $\Phi \in \bcal{K}_{p}$, there exists $\varphi_{\Phi}\in \Pi_{s}$ such that
\begin{eqnarray*}
\mathcal{J}_{\nu}(h,\varphi_{\Phi}) \geq \eta^{\Phi}(h),
\end{eqnarray*}
for any $h\in  \bscr{M}(\KK)$ satisfying 
$\ds  \sup\big\{\eta^{\Phi}(h^{+}) : \Phi\in\bcal{K}_{p} \big\}<+\infty$.
\end{theorem}
\textbf{Proof:} 
For $h\in  \bscr{M}(\KK)$ satisfying  $\ds  \sup\big\{\eta^{\Phi}(h^{+}) : \Phi\in\bcal{K}_{p} \big\}<+\infty$, let us prove the result by showing that
\begin{eqnarray}
\mu^{\varphi_{\Phi}}(h) \geq \eta^{\Phi}(h)
\label{result}
\end{eqnarray}
where $\varphi_{\Phi}$ is the stationary randomized policy induced by $\Phi$ (see \eqref{Def-varphi-Phi}).
There is no loss of generality to assume that $\eta^{\Phi}(h)>-\infty$ and so we have $\eta^{\Phi}(|h|)<\infty$.
We are going to proceed by contradiction to get \eqref{result}.
More precisely, if $\mu^{\varphi_{\Phi}}(h) < \eta^{\Phi}(h)$ then we will introduce a sequence $\{\Psi_{k}\}_{k\in \NN}$ in $\bcal{K}_{p}$ satisfying
$\ds \lim_{k\rightarrow \infty} \eta^{\Psi_{k}}(h)=+\infty$ contradicting the hypothesis. The proof is divided into two steps. We will first introduce $\{\Psi_{k}\}_{k\in \NN}$ and show that
$\Psi_{k}\in \bcal{K}_{p}$ for any $k\in \NN$. In a second step, it will be proven that $\ds \lim_{k\rightarrow \infty} \eta^{\Psi_{k}}(h)=+\infty$ showing the result.

\noindent
\underline{First step: construction of a sequence $\{\Psi_{k}\}_{k\in \NN}$ in $\bcal{K}_{p}$.}
\nl
Let $\mu^{\varphi_{\Phi}}$ be the occupation measure induced by the stationary randomized policy $\varphi_{\Phi}$.
As in the proof of Lemma \ref{strategy-stationary-randomized}, there exists a measurable $\widebar{\RR}_{+}$-valued function $\mathcal{D}_{\varphi_{\Phi}}$ defined on $\XX$ satisfying
\begin{eqnarray}
\mu^{\varphi_{\Phi}}(dx,da)=\mathcal{D}_{\varphi_{\Phi}}(x)  \varphi_{\Phi}(da|x) p(dx).
\label{Def-mu-varphi}
\end{eqnarray}
For $k\in \NN$, consider $\Psi_{k}=(\psi^{\infty},\psi^{*}_{k})$ where $\psi^{\infty}\in \bcal{K}(\AA| \XX)$ is given by
\begin{eqnarray}
\psi^{\infty}(da|x) =  I_{\bcal{E}_{\Phi}^{c}}(x) \varphi_{\Phi}(da|x)
\label{Def-psi-infty}
\end{eqnarray}
and $\psi^{*}_{k}$ is a signed kernel on $\AA$ given $\XX$ defined by
\begin{eqnarray}
\psi^{*}_{k}(da|x) = I_{\bcal{E}_{\Phi}}(x) \Big[ \varphi^{*}(\AA |x) + k \big[ \varphi^{*}(\AA|x)-\mathcal{D}_{\varphi_{\Phi}}(x) \big] \Big]  \varphi_{\Phi}(da|x)
+ (k+1) I_{\bcal{E}_{\Phi}^{c}}(x)  \varphi^{*}(da|x).
\label{Def-psi-*-k}
\end{eqnarray}
Observe that in the previous definition, $\varphi^{*}(\AA|x)-\mathcal{D}_{\varphi_{\Phi}}(x)$ is well defined since $\varphi^{*}\in \bcal{K}(\AA| \XX)$.
To get the result, we will proceed in two steps.
First we will show that $\varphi^{*}(\AA|\cdot)\geq \mathcal{D}_{\varphi_{\Phi}}(\cdot)$ on $\bcal{E}_{\Phi}$ implying that $\psi^{*}_{k}\in \bcal{K}(\AA| \XX)$
and so, $\Psi_{k}\in\bcal{K}(\AA | \XX)^{2}$ for any $k\in \NN$. In a second step, we will prove that $\Psi_{k}\in\bcal{K}_{p}$.

\bigskip

\noindent
$\bullet$ Let us show that $\Psi_{k}\in\bcal{K}(\AA | \XX)^{2}$.
\nl
From \eqref{Def-eta-Phi},  $\eta^{\Phi}(dx,da)= \mathcal{I}_{\infty}(x) \varphi^{\infty}(da|x) p(dx) + \varphi^{*}(da|x) p(dx)$ and so, by using \eqref{Def-E-Phi}
\begin{eqnarray*}
\eta^{\Phi}(dx,da) & = & I_{\bcal{E}_{\Phi}}(x) \varphi^{*}(da|x) p(dx) + I_{\bcal{E}_{\Phi}^{c}}(x) \mathcal{I}_{\infty} (x) \varphi^{\infty}(da |x) p(dx) 
+ I_{\bcal{E}_{\Phi}^{c}}(x)  \varphi^{*}(da|x) p(dx),
\end{eqnarray*}
where by convention $0\times \infty=0$.
Recalling the Definition of $\varphi_{\Phi}$ (see equation \eqref{Def-varphi-Phi}), we easily obtain 
$I_{\bcal{E}_{\Phi}}(x) \varphi^{*}(da|x)=I_{\bcal{E}_{\Phi}}(x) \varphi^{*}(\AA |x) \varphi_{\Phi}(da|x)$ and
$I_{\bcal{E}_{\Phi}^{c}}(x) \mathcal{I}_{\infty} (x) \varphi^{\infty}(da |x)=I_{\bcal{E}_{\Phi}^{c}}(x) \mathcal{I}_{\infty} (x) \varphi_{\Phi}(da|x)$ and so, we get
\begin{align}
\eta^{\Phi}(dx,da) = & I_{\bcal{E}_{\Phi}}(x) \varphi^{*}(\AA |x) \varphi_{\Phi}(da|x) p(dx) + I_{\bcal{E}_{\Phi}^{c}}(x) \mathcal{I}_{\infty} (x)  \varphi_{\Phi}(da|x) p(dx) \nonumber \\
& + I_{\bcal{E}_{\Phi}^{c}}(x)  \varphi^{*}(da|x) p(dx).
\label{mu-Phi-explicit}
\end{align}
Therefore,
\begin{align}
\eta^{\Phi}_{\XX}(dx) & =  \Big[ I_{\bcal{E}_{\Phi}}(x) \varphi^{*}(\AA|x) + I_{\bcal{E}_{\Phi}^{c}}(x) \big[\mathcal{I}_{\infty}(x)+\varphi^{*}(\AA |x)\big]  \Big] p(dx) \nonumber \\
& =  \Big[ I_{\bcal{E}_{\Phi}}(x) \varphi^{*}(\AA|x) + I_{\bcal{E}_{\Phi}^{c}}(x) \mathcal{I}_{\infty}(x) \Big] p(dx).
\label{mu-PhiX-explicit}
\end{align}
Since $\eta^{\Phi}_{\XX}=\nu+\eta^{\Phi}Q$, we have by using \eqref{mu-Phi-explicit}
\begin{align}
\eta^{\Phi}_{\XX}(\Lambda)  = & \nu(\Lambda) +  \int_{\bcal{E}_{\Phi}} Q^{\varphi_{\Phi}}(\Lambda |x) \varphi^{*}(\AA |x) p(dx)
+\int_{\bcal{E}_{\Phi}^{c}} Q^{\varphi_{\Phi}}(\Lambda |x) \mathcal{I}_{\infty}(x) p(dx)  \nonumber \\
& +\int_{\bcal{E}_{\Phi}^{c}} Q^{\varphi^{*}}(\Lambda|x) p(dx),
\label{eq-ref-mu-Phi1}
\end{align}
and with \eqref{mu-PhiX-explicit} it follows
\begin{align*}
\eta^{\Phi}_{\XX}(\Lambda) & = \nu(\Lambda) +  \eta^{\Phi}_{\XX}Q^{\varphi_{\Phi}}(\Lambda) +\int_{\bcal{E}_{\Phi}^{c}} Q^{\varphi^{*}}(\Lambda|x) p(dx).
%\label{eq-ref-mu-Phi2}
\end{align*}
However, $\mu^{\varphi_{\Phi}}_{\XX}$ is the minimal solution to the equation $\beta=\nu+\beta Q^{\varphi_{\Phi}}$ and so, $\mu^{\varphi_{\Phi}}_{\XX}\leq \eta^{\Phi}_{\XX}$.
Combining equations \eqref{Def-mu-varphi} and \eqref{mu-PhiX-explicit}, we obtain $\Big[ I_{\bcal{E}_{\Phi}}(\cdot) \varphi^{*}(\AA|\cdot) + I_{\bcal{E}_{\Phi}^{c}}(\cdot) \mathcal{I}_{\infty}(\cdot) \Big]\geq \mathcal{D}_{\varphi_{\Phi}}(\cdot)$ $p-a.s.$. Consequently, $\mathcal{D}_{\varphi_{\Phi}}(\cdot)\leq \varphi^{*}(\AA|\cdot)$  $p-a.s.$ on $\bcal{E}_{\Phi}$
and according to the definition of $\mathcal{D}_{\varphi_{\Phi}}(\cdot)$ (see equation \eqref{Def-mu-varphi}), there is no loss of generality to claim
\begin{eqnarray}
\mathcal{D}_{\varphi_{\Phi}}(\cdot)\leq \varphi^{*}(\AA|\cdot) \text{ on } \bcal{E}_{\Phi}.
\label{ineq-density}
\end{eqnarray}
Therefore, $\psi^{*}_{k}\in \bcal{K}(\AA| \XX)$ and so, $\Psi_{k}\in\bcal{K}(\AA | \XX)^{2}$ for any $k\in \NN$.

\bigskip

\noindent
\noindent
$\bullet$ Let us show that $\Psi_{k}\in\bcal{K}_{p}$.
\nl
Recalling the definition $\Psi_{k}$ (see equations \eqref{Def-psi-infty}-\eqref{Def-psi-*-k}), we have $\psi^{\infty}(\AA(x)^{c}|x)+\psi^{*}_{k}(\AA(x)^{c}|x)=0$ and 
$\psi^{\infty}(\AA|x)+\psi^{*}_{k}(\AA|x)\geq I_{\bcal{E}_{\Phi}}(x) \varphi^{*}(\AA|x) + I_{\bcal{E}_{\Phi}^{c}}(x) >0$ for any $x\in \XX$.
The only point which remains to prove is that $\eta^{\Psi_{k}}(dx,da)=\mathcal{I}_{\infty}(x) \psi^{\infty}(da|x) p(dx) + \psi^{*}_{k}(da|x) p(dx)$ satisfies
\begin{eqnarray}
\eta^{\Psi_{k}}_{\XX} & = & \nu +\eta^{\Psi_{k}}Q.
\label{eq-ref-gamma}
\end{eqnarray}
Combining the definition of $\Psi_{k}$ (see equations \eqref{Def-psi-infty}-\eqref{Def-psi-*-k}) and the expression of  $\eta^{\Phi}$ (see equation \eqref{mu-Phi-explicit}), we obtain
\begin{eqnarray}
\eta^{\Psi_{k}} & = & \eta^{\Phi}+k\gamma
\label{eta-psi-gamma}
\end{eqnarray}
where $\gamma\in\bcal{M}(\XX\times \AA)$ is given by
\begin{eqnarray}
\gamma(dx,da) = I_{\bcal{E}_{\Phi}}(x) \big[ \varphi^{*}(\AA|x)-\mathcal{D}_{\varphi_{\Phi}}(x) \big]  \varphi_{\Phi}(da|x) p(dx)
+I_{\bcal{E}_{\Phi}^{c}}(x)  \varphi^{*}(da|x) p(dx).
\label{Def-gamma}
\end{eqnarray}

\noindent
To show that \eqref{eq-ref-gamma} holds, we will consider two cases.
\nl
a) Firstly, we will show that equation \eqref{eq-ref-gamma} is satisfied on $\bfrak{B}(\bcal{E}_{\Phi})$. For that, let us consider $\Lambda\in \bfrak{B}(\bcal{E}_{\Phi})$.
From \eqref{eta-psi-gamma}, we have $\eta^{\Psi_{k}}_{\XX}(\Lambda)  = \eta^{\Phi}_{\XX}(\Lambda)+k\gamma_{\XX}(\Lambda) $.
However, $\eta^{\Phi}_{\XX}(\Lambda) = \nu(\Lambda) + \eta^{\Phi} Q(\Lambda)$ showing that $\eta^{\Psi_{k}}_{\XX}(\Lambda) = \nu(\Lambda) + \eta^{\Phi} Q(\Lambda) + k\gamma_{\XX}(\Lambda)$.
If we show that $\gamma_{\XX}(\Lambda) = \gamma Q (\Lambda)$
then $\eta^{\Psi_{k}}_{\XX}(\Lambda) =\nu(\Lambda) +\eta^{\Psi_{k}}Q(\Lambda)$
implying that \eqref{eq-ref-gamma} holds on $\bfrak{B}(\bcal{E}_{\Phi})$.
To see that $\gamma_{\XX}(\Lambda) = \gamma Q (\Lambda)$, observe from \eqref{Def-gamma} that
\begin{eqnarray*}
\gamma_{\XX}(\Lambda) & = & \int_{\Lambda} \big[ \varphi^{*}(\AA|x)-\mathcal{D}_{{\varphi}_{\Phi}}(x) \big]p(dx).
\end{eqnarray*}
Assuming that $\eta^{\Phi}_{\XX}(\Lambda)<\infty$ and combining \eqref{Def-mu-varphi}, \eqref{mu-PhiX-explicit} and the previous equation we have
\begin{eqnarray}
\gamma_{\XX}(\Lambda) & = & \int_{\Lambda} \big[ \varphi^{*}(\AA|x)-\mathcal{D}_{{\varphi}_{\Phi}}(x) \big]p(dx) = \eta^{\Phi}_{\XX}(\Lambda) - \mu^{\varphi_{\Phi}}_{\XX}(\Lambda).
\label{gamma-X=Qgamma-1}
\end{eqnarray}
Now, we obtain by using \eqref{eq-ref-mu-Phi1} and the fact that $\eta^{\Phi}_{\XX}(\Lambda)<\infty$
\begin{eqnarray}
\int_{\bcal{E}_{\Phi}^{c}} Q^{\varphi_{\Phi}}(\Lambda |x) \mathcal{I}_{\infty}(x) p(dx) = 0
\label{Sigma-finite1}
\end{eqnarray}
implying also
\begin{eqnarray}
\int_{\bcal{E}_{\Phi}^{c}} Q^{\varphi_{\Phi}}(\Lambda |x)\mathcal{D}_{\varphi_{\Phi}}(x) p(dx) = 0.
\label{Sigma-finite2}
\end{eqnarray}
Now, combining \eqref{eq-ref-mu-Phi1} and \eqref{Sigma-finite1}
\begin{align*}
\eta^{\Phi}_{\XX}(\Lambda)  & = \nu(\Lambda) +  \int_{\bcal{E}_{\Phi}} Q^{\varphi_{\Phi}}(\Lambda |x) \varphi^{*}(\AA|x) p(dx)
+\int_{\bcal{E}_{\Phi}^{c}} Q^{\varphi^{*}}(\Lambda|x) p(dx).
\end{align*}
Recalling that $\mu^{\varphi_{\Phi}}_{\XX}=\nu+\mu^{\varphi_{\Phi}}_{\XX} Q^{\varphi_{\Phi}}$, we have with \eqref{Def-mu-varphi} and \eqref{Sigma-finite2}
\begin{eqnarray*}
\mu^{\varphi_{\Phi}}_{\XX}(\Lambda) & = & \nu(\Lambda) +  \int_{\bcal{E}_{\Phi}} Q^{\varphi_{\Phi}}(\Lambda |x) \mathcal{D}_{\varphi_{\Phi}}(x) p(dx).
\end{eqnarray*}
The two previous equations gives
\begin{eqnarray}
\eta^{\Phi}_{\XX}(\Lambda) - \mu^{\varphi_{\Phi}}_{\XX}(\Lambda) & = & \int_{\bcal{E}_{\Phi}} Q^{\varphi_{\Phi}}(\Lambda |x) \big[ \varphi^{*}(\AA|x) - \mathcal{D}_{\varphi_{\Phi}}(x) \big]p(dx)
+\int_{\bcal{E}_{\Phi}^{c}} Q^{\varphi^{*}}(\Lambda|x) p(dx).
\label{gamma-X=Qgamma-2}
\end{eqnarray}
From \eqref{gamma-X=Qgamma-1} and \eqref{gamma-X=Qgamma-2}
\begin{eqnarray*}
\gamma_{\XX}(\Lambda) & = & \int_{\bcal{E}_{\Phi}} Q^{\varphi_{\Phi}}(\Lambda |x) \big[ \varphi^{*}(\AA|x) - \mathcal{D}_{\varphi_{\Phi}}(x) \big]p(dx)
+\int_{\bcal{E}_{\Phi}^{c}} Q^{\varphi^{*}}(\Lambda|x) p(dx)
\end{eqnarray*}
Recalling the definition of $\gamma$ (see \eqref{Def-gamma}) we get $\gamma_{\XX}(\Lambda)=\gamma Q(\Lambda)$ for
$\Lambda\in \bfrak{B}(\bcal{E}_{\Phi})$ with $\eta^{\Phi}_{\XX}(\Lambda)<\infty$.
However, equation \eqref{mu-PhiX-explicit} implies that $\eta^{\Phi}_{\XX}$ is $\sigma$-finite on $\bcal{E}_{\Phi}$ and combining
\eqref{mu-Phi-explicit} and \eqref{Def-gamma}, we have $\gamma_{\XX}\leq \eta^{\Phi}_{\XX}$.
Therefore, it follows that $\gamma_{\XX}(\Lambda)=\gamma Q(\Lambda)$ for any $\Lambda$ in $\bfrak{B}(\bcal{E}_{\Phi})$, and so \eqref{eq-ref-gamma} holds on $\bfrak{B}(\bcal{E}_{\Phi})$.

\bigskip

\noindent
b) Secondly, we will show that equation \eqref{eq-ref-gamma} is satisfied on $\bfrak{B}(\bcal{E}_{\Phi}^{c})$. For that, let $\Lambda\in \bfrak{B}(\bcal{E}_{\Phi}^{c})$. It is important to observe from \eqref{mu-PhiX-explicit} that in this case $\eta^{\Phi}_{\XX}(\Lambda)=0$ or $+\infty$.
Therefore, we obtain on one hand 
$\eta^{\Psi_{k}}_{\XX}(\Lambda)= \eta^{\Phi}_{\XX}(\Lambda)+k\gamma_{\XX}(\Lambda)=\eta^{\Phi}_{\XX}(\Lambda)$ by recalling 
\eqref{eta-psi-gamma} and using the fact that $\gamma_{\XX}\leq \eta^{\Phi}_{\XX}$
and on the other hand $\eta^{\Phi}_{\XX}(\Lambda)=\eta^{\Phi}_{\XX}(\Lambda)+k\gamma Q(\Lambda)$ since by \eqref{Def-gamma}
\begin{eqnarray*}
\gamma Q (\Lambda)  & = & \int_{\bcal{E}_{\Phi}} Q^{\varphi_{\Phi}}(\Lambda |y) \big[ \varphi^{*}(\AA|y)-\mathcal{D}_{\varphi_{\Phi}}(y) \big] p(dy)
+\int_{\bcal{E}_{\Phi}^{c}} Q^{\varphi^{*}}(\Lambda|y) p(dy) \leq \eta^{\Phi}_{\XX}(\Lambda)
\end{eqnarray*}
where the last inequality comes from \eqref{eq-ref-mu-Phi1}.
Therefore,
$\eta^{\Psi_{k}}_{\XX}(\Lambda)=\eta^{\Phi}_{\XX}(\Lambda)+k\gamma Q(\Lambda)=\nu\Lambda)+\eta^{\Psi_{k}}Q(\Lambda)$
showing that \eqref{eq-ref-gamma} holds on $\bfrak{B}(\bcal{E}_{\Phi}^{c})$.

\bigskip

Finally, equation \eqref{eq-ref-gamma} is satisfied and as a consequence $\Psi_{k} \in \bcal{K}_{p}$ for any $k\in \NN$.

\bigskip

\noindent
\underline{Second step: $\ds \lim_{k\rightarrow \infty} \eta^{\Psi_{k}}(h)=+\infty$.}
\nl
Recalling that $\eta^{\Phi}(|h|)<\infty$, we get from \eqref{mu-Phi-explicit} 
\begin{eqnarray*}
\int_{\bcal{E}_{\Phi}^{c}} |h(x,a)| \mathcal{I}_{\infty} (x) \varphi_{\Phi}(da |x) p(dx) =0
\end{eqnarray*}
implying also
\begin{eqnarray*}
\int_{\bcal{E}_{\Phi}^{c}} |h(x,a)| \mathcal{D}_{\varphi_{\Phi}}(x) \varphi_{\Phi}(da |x) p(dx) =0.
\end{eqnarray*}
Therefore, combining \eqref{Def-mu-varphi}, \eqref{mu-Phi-explicit}, \eqref{Def-gamma} and the two previous equations we obtain easily that
\begin{eqnarray*}
\eta^{\Phi} (h) - \mu^{\varphi_{\Phi}}(h) = \gamma(h).
\end{eqnarray*}
If $\eta^{\Phi}(h)> \mu^{\varphi_{\Phi}}(h)$ then $\gamma(h)>0$ and $\ds \lim_{k\rightarrow \infty}\eta^{\Psi_{k}}(h)=\eta^{\Phi}(h)+  \lim_{k\rightarrow \infty}k\gamma(h)=+\infty$
giving the result.
\hfill$\Box$

%%%%%%%%%%%%%%%
\subsection{Properties of the constrained control problem}
\label{sec-sufficiency}
%%%%%%%%%%%%%%%
The main objective of this subsection is to show that in the framework of constrained control problems, the supremum of the expected total rewards over the set of randomized policies is equal to the supremum of the expected total rewards over the set of \textit{stationary} randomized policies.
Our results use Theorem \ref{Schal-Theorem} presented in the Appendix which is a slight modification of Theorem 1 in  Sch\"{a}l \cite{schal83} who has established a stronger version of this type of result but in the unconstrained case. To use Sch\"{a}l's results, we need to impose Conditions (W) or (S) and in addition, to deal with the constrained case, we need to impose a Slater-type condition.

\bigskip

The next technical Lemma shows that, roughly speaking, under Assumption \ref{Hyp-finiteness-linear}, the unconstrained control problems associated to a reward function given by either $r$ or $c_{i}$ for $i\in\NN_{q}$ are different from $+\infty$.
\begin{lemma}
\label{Finiteness-occupation}
Suppose Assumptions \ref{Hyp-Transition-kernel-absolute-continuity} and \ref{Hyp-finiteness-linear} and either Conditions (W) or (S) hold. Then,
$$\sup \big\{ \mu (r^{+}) : \mu\in\bcal{O} \big\} < +\infty \text{ and } \sup \big\{\mu (c^{+}_{i}) : \mu\in\bcal{O} \big\}<+\infty,$$
for $i\in\NN_{q}$.
\end{lemma}
\textbf{Proof:} 
The idea is to apply Theorem \ref{Schal-Theorem} to the unconstrained models associated to the reward functions given by one of the following mappings:
$r^{+}$ and $c^{+}_{i}$ for $i\in\NN_{q}$.
Clearly, the Convergence Assumption and the Continuity and Compactness Assumptions in \cite[p. 367]{schal83} are satisfied. Therefore, we have by using Theorem \ref{Schal-Theorem}
\begin{align*}
\sup\{ \mu(h) :\mu\in\bcal{O} \} =  \sup\{ \mu(h) :\mu\in\bcal{O}_{s} \}
\end{align*}
for any function $h$ given by either $r^{+}$ or $c^{+}_{i}$ for $i\in\NN_{q}$. Now, from Assumption \ref{Hyp-Transition-kernel-absolute-continuity} we can apply
Lemma \ref{strategy-stationary-randomized} to have
\begin{align*}
\sup\{ \mu(h) :\mu\in\bcal{O}_{s} \} \leq \sup\{ \eta^{\Phi} (h) : \Phi\in\bcal{K}_{p}\}.
\end{align*}
Recalling Assumption \ref{Hyp-finiteness-linear} we obtain the result.
\hfill$\Box$

\bigskip

The next result shows that if the Slater condition is satisfied for an arbitrary policy then there exists a stationary randomized policy satisfying the same type of condition.
\begin{proposition}
\label{Slater-general=stationary}
Suppose Assumptions \ref{Hyp-Transition-kernel-absolute-continuity}, \ref{Hyp-finiteness} and either Conditions (W) or (S) hold.
If the Slater condition is satisfied, then there exists $\widetilde{\mu}\in\bcal{O}_{s}$ satisfying $\theta^{*}_{i}<\widetilde{\mu}(c_{i})$
for any $i\in\NN_{q}$.
\end{proposition}
\textbf{Proof:} 
The result is proved by induction.
Applying Theorem \ref{Schal-Theorem} for the unconstrained model associated to the reward function $c_{1}$, we have
\begin{eqnarray*}
\sup\{ \mu (c_{1}) : \mu\in \bcal{O} \}  = \sup\{ \mu (c_{1}) : \mu\in \bcal{O}_{s} \}.
\end{eqnarray*}
Since $\mu^{*}(c_{1})>\theta_{1}^{*}$ (by recalling the Slater condition), we have $\sup\{ \mu (c_{1}) : \mu\in \bcal{O}_{s} \}>\theta_{1}^{*}$ implying the existence of 
$\mu_{1}\in \bcal{O}_{s}$ such that $\mu_{1}(c_{1})>\theta_{1}^{*}$.
For $n\in \NN_{q-1}$, let us assume the existence of $\mu_{n}\in \bcal{O}_{s}$ such that $\mu_{n}(c_{i})>\theta_{i}^{*}$ for $i\in \NN_{n}$.
Therefore, we can combine Lemma \ref{Finiteness-occupation} and Proposition \ref{constraint=unconstraint-Slater-stationary} to obtain
\begin{align*}
\sup\big\{ \mu (c_{n+1}) : \mu\in \bcal{O} \text{ and } \mu(c_{i})>\theta_{i}^{*} & \text{ for } i\in \NN_{n} \big\} \nonumber \\
& = \sup\big\{ \mu (c_{n+1}) : \mu\in \bcal{O}_{s} \text{ and } \mu(c_{i})>\theta_{i}^{*} \text{ for } i\in \NN_{n} \big\}.
\end{align*}
However,
$$\sup\Big\{ \mu (c_{n+1}) : \mu\in \bcal{O} \text{ and }\mu(c_{i})>\theta_{i}^{*} \text{ for } i\in \NN_{n}  \Big\} \geq \mu^{*}(c_{n+1}) > \theta_{n+1}^{*}$$
implying the existence of $\mu_{n+1}\in \bcal{O}_{s}$ such that $\mu_{n+1}(c_{i})>\theta_{i}^{*}$ for $i\in \NN_{n+1}$. This gives the result.
\hfill$\Box$

\bigskip

Below is the main result of this subsection that states roughly speaking that in the framework of constrained control problems, the supremums of the expected total rewards over the set of randomized policies and over the set of stationary randomized policies coincide.
\begin{theorem}\label{theorem-constraint=unconstraint}
Suppose Assumptions \ref{Hyp-Transition-kernel-absolute-continuity}, \ref{Hyp-finiteness} and either Conditions (W) or (S) hold.
If the Slater condition is satisfied, then 
\begin{align*}
\sup\{ \mu(r) :\mu\in\bcal{O} \hbox{ and } \mu(c_{i})\geq\theta^{*}_{i} & \text{ for } i\in\NN_{q}\} \nonumber \\
& =  \sup\{ \mu(r) :\mu\in\bcal{O}_{s} \hbox{ and } \mu(c_{i})\geq\theta^{*}_{i} \text{ for } i\in\NN_{q}\}.
\end{align*}
\end{theorem}
\textbf{Proof:}
Applying Proposition \ref{Slater-general=stationary}, there exists $\widetilde{\mu}\in\bcal{O}_{s}$ satisfying the Slater condition, that is, $\widetilde{\mu}(c_{i})>\theta^{*}_{i}$
for $i\in\NN_{q}$.
Now, combining Lemma \ref{Finiteness-occupation} and Proposition \ref{constraint=unconstraint-Slater-stationary}, we obtain the result.
\hfill$\Box$

%%%%%%%%%%%%%%%%%%%%%%%%%%%%%%%
\section{Main results}
\label{LP-formulation=Control-Problem}
%%%%%%%%%%%%%%%%%%%%%%%%%%%%%%%
In this section, we present the main results of this paper showing that the original control problem is equivalent to the convex program introduced in Definition \ref{Def-Linear-Program} for a weakly or strongly continuous transition kernel.

%%%%%%%%%%%%%%%%%%%%%%%%%%%%%%%
\paragraph{The case of Condition (W)}
%%%%%%%%%%%%%%%%%%%%%%%%%%%%%%%
\begin{theorem}
\label{theorem-main-LP-condition-W}
Suppose Assumptions \ref{Hyp-Transition-kernel-absolute-continuity}, \ref{Hyp-finiteness} and Condition (W) hold.
If the Slater condition is satisfied, then
\begin{align}
\sup\big\{  \mathcal{J}_{\nu}(r,\pi) :\pi\in \Pi_{\theta^{*}} \big\}
& = \sup\big\{ \eta^{\Phi}(r) : \Phi\in\bcal{K}_{p} \text{ and } \eta^{\Phi}(c_{i}) \geq \theta^{*}_{i} \text{ for } i\in\NN_{q} \big\}
\label{Values=condition-W}
\end{align}
where $p\in\bcal{P}(\XX)$ is defined in \eqref{Def-p}.
Moreover,  if $\hat{\Phi}$ is an optimal solution to the convex program $\bcal{KP}_{p}$ then
the stationary randomized policy $\varphi_{\hat{\Phi}}$ induced by $\hat{\Phi}$ is optimal
for the constrained control problem, that is,
\begin{align}
\mathcal{J}_{\nu}(r,\varphi_{\hat{\Phi}})=\sup\big\{  \mathcal{J}_{\nu}(r,\pi) :\pi\in \Pi_{\theta^{*}} \big\}.
\end{align}
\end{theorem}
\textbf{Proof:} Theorem \ref{theorem-constraint=unconstraint} states that
\begin{align*}
\sup\big\{  \mathcal{J}_{\nu}(r,\pi) :\pi\in \Pi_{\theta^{*}} \big\}
& = \sup\big\{  \mathcal{J}_{\nu}(r,\pi) : \pi\in\Pi_{s} \inter  \Pi_{\theta^{*}} \big\}.
\end{align*}
However, from Lemma \ref{strategy-stationary-randomized}, we have 
\begin{align*}
\sup\big\{ \mathcal{J}_{\nu}(r,\pi) :  \pi\in\Pi_{s} \inter  \Pi_{\theta^{*}} \big\}  
& \leq \sup\big\{ \eta^{\Phi}(r) : \Phi\in\bcal{K}_{p} \text{ and } \eta^{\Phi}(c_{i}) \geq \theta^{*}_{i} \text{ for } i\in\NN_{q} \big\},
\end{align*}
Now, consider $\Phi \in \bcal{K}_{p}$. By using Theorem \ref{Linear-prog=stationary-optimality}, $\mathcal{J}_{\nu}(h,\varphi_{\Phi}) \geq \eta^{\Phi}(h)$ for $h$ given either $r$ or $c_{i}$ for $i\in\NN_{q}$ implying that
$\varphi_{\Phi}\in \Pi_{s}\inter  \Pi_{\theta^{*}}$ and also the reverse inequality
\begin{align*}
\sup\big\{ \mathcal{J}_{\nu}(r,\pi) :  \pi\in\Pi_{s} \inter  \Pi_{\theta^{*}} \big\}  
& \geq \sup\big\{ \eta^{\Phi}(r) : \Phi\in\bcal{K}_{p} \text{ and } \eta^{\Phi}(c_{i}) \geq \theta^{*}_{i} \text{ for } i\in\NN_{q} \big\}
\end{align*}
showing the first part of the result.

Now if $\hat{\Phi} \in\bcal{K}_{p}$ is an optimal solution to the convex program $\bcal{KP}_{p}$ then
$\eta^{\hat{\Phi}}(c_{i}) \geq \theta^{*}_{i}$ for any $i\in\NN_{q}$ and
$\eta^{\hat{\Phi}}(r) = \sup\big\{ \eta^{\Phi}(r) : \Phi\in\bcal{K}_{p} \text{ and } \eta^{\Phi}(c_{i}) \geq \theta^{*}_{i} \text{ for } i\in\NN_{q} \big\}$.
Therefore, the stationary randomized policy
$\varphi_{\hat{\Phi}}\in \Pi_{\theta^{*}}$ satisfies $\mathcal{J}_{\nu}(r,\varphi_{\hat{\Phi}}) \geq \eta^{\hat{\Phi}}(r)$ by using Theorem \ref{Linear-prog=stationary-optimality}.
Now, by using the first part of the result (see equation \eqref{Values=condition-W}) it follows that
$\mathcal{J}_{\nu}(r,\varphi_{\hat{\Phi}}) \geq \sup\big\{ \mathcal{J}_{\nu}(r,\pi) :  \pi\in \Pi_{\theta^{*}} \big\} $
giving the last part of the result.
\hfill$\Box$

\begin{remark}
As mentioned in the introduction, the previous result has the advantage of proposing a convex programming formulation for constrained MDPs under the ETR criterion with signed reward functions and satisfying condition (W) which has not been so far addressed in the literature.
In \cite{dufour12}, the authors do not really analyse a convex program, but study a related optimization problem where the MPDs under consideration satisfy condition (W) but the proposed approach strongly relies on the positiveness of the cost functions and cannot be generalized to the framework of signed cost functions.
\end{remark}
%%%%%%%%%%%%%%%%%%%%%%%%%%%%%%%
\paragraph{The case of condition (S)}
%%%%%%%%%%%%%%%%%%%%%%%%%%%%%%%

\begin{theorem}\label{theorem-main-LP-condition-S}
Suppose Assumptions \ref{Hyp-finiteness} and Condition (S) hold.
If the Slater condition is satisfied, then
\begin{align*}
\sup\big\{  \mathcal{J}_{\nu}(r,\pi) :\pi\in \Pi_{\theta^{*}} \big\}
& = \sup\big\{ \eta^{\Phi}(r) : \Phi\in\bcal{K}_{p} \text{ and } \eta^{\Phi}(c_{i}) \geq \theta^{*}_{i} \text{ for } i\in\NN_{q} \big\}
\end{align*}
where $p$ is defined in \eqref{Def-p} for $P$ given by \eqref{Def-P-strongly-continuous}.
Moreover,  if $\hat{\Phi}$ is an optimal solution to the convex program 
then the stationary randomized policy $\varphi_{\hat{\Phi}}$ induced by $\hat{\Phi}$ is optimal
for the constrained control problem introduced in Definition \ref{Constrained-control}.
\end{theorem}
\textbf{Proof:}
{Up to the definition of $p$ whose existence is established in Lemma \ref{Existence-P-strongly-continuous},} the proof
of this result is identical to that of Theorem \ref{theorem-main-LP-condition-W}.
\hfill$\Box$

\begin{remark}
In \cite{dufour13}, the authors  do not really analyse a convex program but study a related optimization problem where the MPDs under consideration satisfy condition (S) by assuming that the transition kernel is absolutely continuous with respect to a reference probability measure uniformly in the state and action variables. In the previous result, we show that this assumption is not needed under condition (S) if this hypothesis is replaced by a Slater-type condition.
\end{remark}

%%%%%%%%%%%%%%%%%%%%%%%%%%%%%%%
\section{Example}
\label{sec-example}
%%%%%%%%%%%%%%%%%%%%%%%%%%%%%%%
In this section, we provide an example with one constraint to illustrate our results and compare them with reference \cite{dufour13}.
The results obtained in \cite{dufour12} cannot be used for this model because the contraint function takes positive and negative values.
%Actually, this example is borrowed from \cite{dufour13} (see section 6.2).
We will show that one of the conditions of \cite{dufour13} is not satisfied while the approach developed in the present paper can be applied.
This example shows that there is a gap between the initial optimization problem and  the mathematical program associated to the measures satisfying the characteristic equation, that is,
\begin{align*}
 \sup\big\{ \mathcal{J}_{\nu}(r,\pi) : \pi\in\Pi  \text{ and } & \mathcal{J}_{\nu}(c_{1},\pi) \geq \theta^{*}_{1}  \big\} \\
 & < \sup\big\{ \mu(r) : \mu\in\bcal{M}(\mathbf{X}),~\mu_{\XX}= \nu + \mu Q \text{ and } \mu(c_1)\geq \theta^\ast_1\big\}.
\end{align*}
It means that the characteristic equation $\mu_{\XX}= \nu + \mu Q$  generates measures that do not correspond to any occupation measures of the process. This type of measures has been called in \cite{dufour10} phantom solutions of the characteristic equation.
The interesting point is that at the same time, we may have 
\begin{align*}
\sup\big\{ \mathcal{J}_{\nu}(r,\pi) : \pi\in\Pi  \text{ and }\mathcal{J}_{\nu}(c_{1},\pi) \geq \theta^{*}_{1}  \big\} = \sup\big\{ \eta^{\Phi}(r) : \Phi\in\bcal{K}_{p}\text{ and }\eta^\Phi(c_1)\geq \theta^\ast_1\big\}.
\end{align*}
This means that the set $\big\{ \eta^{\Phi} : \Phi\in\bcal{K}_{p}\big\}$ which is by the way a subset of $\big\{ \mu\in\bcal{M}(\mathbf{X}) : \mu_{\XX}= \nu + \mu Q \}$ may generate less of phantom solutions.

Two different values of the constraint limit $\theta^\ast_1$ will be studied.
For the first value of the constraint limit, it will be shown that the approach proposed in the present paper can be applied implying that the value of the original control problem coincides with the value of the convex program $\bcal{KP}_{p}$.
When changing the value of the constraint limit, the Slater condition will not be satisfied.
However, it is interesting to observe that in this latter case, the values of the original control problem and its associated convex program $\bcal{KP}_{p}$ still coincide although the Slater condition is not fulfilled.
It appears that the Slater condition is not a necessary condition to establish the correspondance between the constrained control problem and its associated
convex program $\bcal{KP}_{p}$.

\bigskip

We consider the control model 
\begin{equation*}
\big(\mathbf{X},\mathbf{A},Q,r,c_1,\theta^{*}_1,\nu\big) 
\end{equation*}
where $\mathbf{X}=\mathbb{Z}\cup\{\Delta\}$ and the action set is given by $\mathbf{A}=\{a,b\}$. 
For $x\neq 1$, $\mathbf{A}(x)=\{a\}$; $\mathbf{A}(1)=\{a,b\}$ and $\mathbf{A}(\Delta)=\{a\}$.
The stochastic kernel $Q$ is given by $Q(x+1|x,a)=1$ for $x\leq 0$ and
$Q(y|x,a)=(1/2) I_{\{x+1\}}(y) + (1/2) I_{\{x+2\}}(y) $, for $x\geq 1$ and finally, $Q(\Delta |1,b)=Q(\Delta |\Delta,a)=1$.
The one-step reward function is given by $r(x,a)=(1/2)^{|x|}$ for $x\neq 1$; $r(1,a)=r(1,b)=1/2$ and $r(\Delta,a)=0$.
The one-step constraint function is given by $c_{1}(x)=(-1/2)^{|x|}$ for $x\neq 1$; $c_{1}(1,a)=-1/18$ and $c_{1}(1,b)=1$.
The initial distribution $\nu$ satisfies $\nu(\{1\})=\nu(\{\Delta\})=1/2$.
The constraint limit is given by $\theta^{*}_{1}$. Two cases are studied: $\theta^{*}_{1}=1/4$ and $\theta^{*}_{1}=1/2$.

\bigskip

Let $\mu\in\bcal{M}(\mathbf{X})$ satisfying the characteristic equation $\mu_{\XX}= \nu + \mu Q$ and so,
$\mu(\Delta,a)=+\infty$; $\mu(x,a)=\mu(0,a)$ for $x\leq 0$; $\mu(1,a)+\mu(1,b)=1/2+\mu(0,a)$ and finally, $\mu(2,a)=(1/2)\mu(1,a)$ and
$\mu(x,a)=(1/2)\mu(x-1,a)+(1/2)\mu(x-2,a)$ for $x\geq 3$ showing that for $x\geq 2$,
$\mu(x,a)=(1/6)[4-(-1/2)^{x-2}] \mu(1,a)$.
Therefore, 
\begin{align*}
\sup \big\{\mu(r) : \mu\in\bcal{M}(\mathbf{X}),~\mu_{\XX}= \nu + \mu Q \big\} \geq  \sup \big\{\mu(0,a) : \mu(0,a) \in \widebar{\RR}_{+} \big\} = +\infty
\end{align*}
since $\ds \mu(r)=\sum_{x\neq 1} (1/2)^{|x|} \mu(x,a)+ (1/2) [\mu(1,a)+\mu(1,b)]$. This implies that Assumption (A2) in \cite{dufour13} is not satisfied and therefore, the approach developed there cannot be applied.

\bigskip

The stochastic kernel $P$ on $\mathbf{X}$ given $\mathbf{X}$ defined by $P(x|y)=Q(x|y,a)$ for $y\in\{\Delta\}\cup\ZZ\setminus\{1\}$
and $P(2|1)=P(3|1)=P(\Delta |1)=1/3$ satisfies Assumption \ref{Hyp-Transition-kernel-absolute-continuity}.

The probability $p$ associated to $P$ and given by \eqref{Def-p} satisfies $p(x)=0$ for $x\leq 0$.
As a consequence, $\eta^\Phi(x,a)=0$ for any $x\leq 0$ and $\Phi\in\bcal{K}_{p}$. Moreover, since $\eta^\Phi$ satisfies the characteristic equation, it follows that
$\eta^\Phi(1,a)+\eta^\Phi(1,b)=1/2$ and $\eta^\Phi(x,a)=(1/6)[4-(-1/2)^{x-2}] \eta^\Phi(1,a)$ for $x\geq 2$
and $\eta^\Phi(\Delta,a)=+\infty$.
Thus, 
\begin{align*}
\eta^{\Phi}(r) & = r(1,a) \eta^\Phi(1,a)+ r(1,b) \eta^\Phi(1,b) + \sum_{x\geq 2} r(x,a) \eta^\Phi(x,a)\\
& = (1/2) [\eta^\Phi(1,a)+\eta^\Phi(1,b)] + \sum_{x\geq 2} (1/2)^{x} (1/6)[4-(-1/2)^{x-2}] \eta^\Phi(1,a) \\
& = 1/4+ (3/10) \: \eta^\Phi(1,a)
\end{align*}
and similarly,
\begin{align*}
\eta^{\Phi}(c_{1}) & = (-1/18) \eta^\Phi(1,a) +\eta^\Phi(1,b) + \sum_{x\geq 2} (-1/2)^{x} (1/6)[4-(-1/2)^{x-2}] \eta^\Phi(1,a) \\
& = 1/2- \eta^\Phi(1,a)
\end{align*}
where $\eta^\Phi(1,a)\in[0,1/2]$.
Clearly, we have $\eta^{\Phi}(r^{+})<+\infty$ and $\eta^{\Phi}(c^{+}_{1})<+\infty$ for any $\Phi\in\bcal{K}_{p}$ showing that Assumption \ref{Hyp-finiteness-linear} is satisfied.

Now, let $\pi_{a}$ (respectively, $\pi_{b}$) be the deterministic stationary policy given by $\pi_{a}(\{a\}|x)=1$ for $x\in \ZZ\union\{\Delta\}$ (respectively, $\pi_{b}(\{a\}|x)=1$ if $x\in \ZZ\union\{\Delta\} \setminus\{1\}$ and $\pi_{b}(\{b\}|1)=1$).
It is easy to see that the occupation measure $\mu^{\pi_{a}}$ is given by $\mu^{\pi_{a}}(1,a)=1/2$; $\mu^{\pi_{a}}(1,b)=0$; $\mu^{\pi_{a}}(\Delta,a)=+\infty$;
$\mu^{\pi_{a}}(x,a)=0$ for any $x\leq 0$ and
$\mu^{\pi_{a}}(x,a)=(1/12)[4-(-1/2)^{x-2}]$ for $x\geq 2$
and the occupation measure $\mu^{\pi_{b}}$ satisfies $\mu^{\pi_{b}}(x,a)=0$ for any $x\in \ZZ$; $\mu^{\pi_{b}}(1,b)=1/2$ and $\mu^{\pi_{b}}(\Delta,a)=+\infty$.
It follows easily $\mu^{\pi_{a}}(r)=\sum_{x\geq 2} (1/2)^{x} (1/12)[4-(-1/2)^{x-2}] + 1/4=2/5$ and $\mu^{\pi_{b}}(r)=r(1,b) \mu^{\pi_{b}}(1,b)=1/4$.
Observe also that $\mu^{\pi_{a}}(c_1)= -1/18  + \sum_{x\geq 2} (-1/2)^{x} (1/6)[4-(-1/2)^{x-2}] =0$ and
$\mu^{\pi_{b}}(c_1)=1/2$.
Clearly, the reward $\mathcal{J}_{\nu}(r,\pi)$ takes values in the interval $[\mathcal{J}_{\nu}(r,\pi_{b}),\mathcal{J}_{\nu}(r,\pi_{a})]$ when the policy $\pi$ ranges over
$\Pi$ and the constraint $\mathcal{J}_{\nu}(c_1,\pi)$ takes values in $[\mathcal{J}_{\nu}(c_{1},\pi_{a}),\mathcal{J}_{\nu}(c_{1},\pi_{b})]$.
Therefore, Assumption \ref{Hyp-negative-costs-finiteness} is satisfied.

Finally, Condition (W) is obviously satisfied for this model.

\bigskip

Remark that for any $\alpha\in[0,1]$, the stationary randomized policy given by $\pi(\{a\}|1)=\alpha$, $\pi(\{b\}|1)=1-\alpha$ and $\pi(\{a\}|x)=1$ for $x\in\ZZ\setminus\{1\}$ yields $\mathcal{J}_{\nu}(r,\pi)=(1-\alpha) \mathcal{J}_{\nu}(r,\pi_{b})+\alpha \mathcal{J}_{\nu}(r,\pi_{a})$
and $\mathcal{J}_{\nu}(c_{1},\pi)=(1-\alpha) \mathcal{J}_{\nu}(c_{1},\pi_{b})+\alpha \mathcal{J}_{\nu}(c_{1},\pi_{a})$.

\bigskip

\noindent\underline{The case where $\theta^{*}_{1}=1/4.$}
From the previous discussion, we have
\begin{align*}
 \sup\big\{ \mathcal{J}_{\nu}(r,\pi) : \pi\in\Pi  \text{ and }\mathcal{J}_{\nu}(c_{1},\pi) \geq \theta^{*}_{1}  \big\} = \mathcal{J}_{\nu}(r,\pi^{*}) = 13/40
\end{align*}
where $\pi^{*}$ is the stationary randomized policy given by $\pi^{*}(\{a\}|1)=\pi^{*}(\{b\}|1)=1/2$, $\pi^{*}(\{a\}|x)=1$ for $x\in\ZZ\setminus\{1\}$.
Moreover,
\begin{align*}
\sup\big\{ \eta^{\Phi}(r) & : \Phi\in\bcal{K}_{p}\text{ and }\eta^\Phi(c_1)\geq \theta^\ast_1\big\} \\
& =\sup\{ (3/10) \: \eta^\Phi(1,a) + 1/4 : \eta^\Phi(1,a) \in [0,1/2] \text{ and } (1/2-\eta^\Phi(1,a) )\geq 1/4\} \\
& = 13/40.
\end{align*}
Therefore, the values of the original control problem and the convex program $\bcal{KP}_{p}$ agree
as claimed by Theorem \ref{theorem-main-LP-condition-W} since the Slater condition holds.

{Observe that the optimal value of the convex program $\bcal{KP}_{p}$ is achieved for $\eta^{\hat{\Phi}}(1,a)=1/4$ where $\hat{\Phi} \in\bcal{K}_{p}$ is an optimal solution to the convex program $\bcal{KP}_{p}$. Since $p(1)=1/4$, the stationary policy $\varphi_{\hat\Phi}$ induced by $\hat\Phi$ is given by $\varphi_{\hat\Phi}(\{a\}|1)=\varphi_{\hat\Phi}(\{b\}|1)=1/2$ and $\varphi_{\hat\Phi}(\{a\} |\Delta)=\varphi_{\hat\Phi}(\{a\}|x)=1$ for $x\in\ZZ\setminus\{1\}$. This optimal policy corresponds to $\pi^{*}$ as determined above.}

\bigskip

\bigskip

\noindent\underline{The case where $\theta^{*}_{1}=1/2.$} We have for this value of the constraint limit,
\begin{align*}
 \sup\big\{ \mathcal{J}_{\nu}(r,\pi) : \pi\in\Pi  \text{ and }\mathcal{J}_{\nu}(c_{1},\pi) \geq \theta^{*}_{1}  \big\} = \mathcal{J}_{\nu}(r,\pi^{*}) = 1/4
\end{align*}
where $\pi^{*}$ is the stationary randomized policy given by $\pi^{*}(\{a\}|1)=0$, $\pi^{*}(\{b\}|1)=1$, $\pi^{*}(\{a\}|x)=1$ for $x\in\ZZ\setminus\{1\}$.

However, we cannot apply the results of the present paper because in this case the Slater condition is not satisfied. Indeed, for any $\pi\in \Pi$,
$\mathcal{J}_{\nu}(c_{1},\pi)\leq 1/2$. But, the values of the original control problem and the convex program $\bcal{KP}_{p}$ still agree since
\begin{align*}
\sup\big\{ \eta^{\Phi}(r) & : \Phi\in\bcal{K}_{p}\text{ and }\eta^\Phi(c_1)\geq \theta^\ast_1\big\} \\
& =\sup\{ (3/10) \: \eta^\Phi(1,a) + 1/4 : \eta^\Phi(1,a) \in [0,1/2] \text{ and } (1/2-\eta^\Phi(1,a) )\geq 1/2\} \\
& = 1/4.
\end{align*}

%%%%%%%%%%%%%%%%%%%%%%%%%%%%%%%
%%%%%%%%%%%%%%%%%%%%%%%%%%%%%%%
\appendix 
%%%%%%%%%%%%%%%%%%%%%%%%%%%%%%%
%%%%%%%%%%%%%%%%%%%%%%%%%%%%%%%

%%%%%%%%%%%%%%%%%%%%%%%%%%%%%%%
\section{Appendix}
%%%%%%%%%%%%%%%%%%%%%%%%%%%%%%%
In this appendix, let $m$ be an integer in $\NN^{*}$. Consider the functions $h\in \bscr{M}(\KK)$ and $g_{i}\in \bscr{M}(\KK)$ for $i\in\NN_{m}$.
We will first present a slightly different version of a result derived by M. Sch\"{a}l in \cite[Theorem 1]{schal83}.
The only difference is that, we consider here the expected total reward  criterion while in \cite{schal83}, Sch\"al deals with the conditional version of that performance criterion.
We will use it repeatedly in this paper. In this section we will also establish a technical result that is used in section \ref{sec-sufficiency} to show that in the framework of control problems with constraints, the supremum of the expected total rewards over the set of randomized policies is equal to the supremum of the expected total rewards over the set of \textit{stationary} randomized policies.

To use Theorem 1 in \cite{schal83}, we need to introduce the following two sets of conditions:
\begin{enumerate}[label=\roman*)]
\item[$\mathbf{(\bcal{S}1)}$] For any $x\in \XX$, $\AA(x)$ is compact.
\item[$\mathbf{(\bcal{S}2)}$] For any $x\in \XX$ and $\Lambda\in \bfrak{B}(\XX)$, $Q(\Lambda |x,\cdot)$ is continuous on $\AA(x)$.
\item[$\mathbf{(\bcal{S}3)}$] For any $x\in \XX$, $h(x,\cdot)$ is upper-semicontinuous on $\AA(x)$.
\item[$\mathbf{(\bcal{S}4)}$] For any $x\in \XX$, $g_{i}(x,\cdot)$ for $i\in\NN_{m}$ are upper-semicontinuous on $\AA(x)$.
\end{enumerate}
or
\begin{enumerate}[label=\roman*')]
\item[$\mathbf{(\bcal{W}1)}$] For any $x\in \XX$, the action set $\AA(x)$ is compact and the multifunction from $\XX$ to~$\AA$ defined by $x\rightarrow \AA(x)$ is upper-semicontinuous.
\item[$\mathbf{(\bcal{W}2)}$] For any $f\in \bcal{C}(\XX)$, $Qf$ is continuous on $\KK$.
\item[$\mathbf{(\bcal{W}3)}$] The function $h$ is upper-semicontinuous on $\KK$.
\item[$\mathbf{(\bcal{W}4)}$] The functions $g_{i}$ for $i\in\NN_{m}$ are upper-semicontinuous on $\KK$.
\end{enumerate}

\begin{theorem}
\label{Schal-Theorem}
Suppose $\mu(h^{+})< +\infty$ or $\mu(h^{-})< +\infty$ for any $\mu\in \bcal{O}$ and either conditions $(\mathcal{S}1)$-$(\mathcal{S}3)$ or $(\mathcal{W}1)$-$(\mathcal{W}3)$ are satisfied.
Then 
\begin{align}
\sup\big\{ \mu(h) :\mu\in\bcal{O}\big\} =  \sup\big\{ \mu(h) :\mu\in\bcal{O}_{s} \big\}.
\end{align}
\end{theorem}
\textbf{Proof:} The proof of this result is essentially the same as Theorem 1 in \cite{schal83}. The only difference is that, we consider here the expected total reward  criterion while in \cite{schal83}, Sch\"al deals with the conditional version of that performance criterion. By adapting the arguments developed in \cite{schal83}, we obtain easily the result.
\hfill$\Box$

\begin{proposition}\label{constraint=unconstraint-Slater-stationary}
Consider $\tilde{\theta}\in \RR^{m}$. 
{Assume $\sup \big\{ \mu(h^{+}+g_{i}^{+}) : \mu\in\bcal{O}\union\{\eta^{\Phi}: \Phi\in\bcal{K}_{p}\} \big\}<+\infty$; 
$\mu(h^{-})< +\infty$ and $\mu( g_{i}^{-})<+\infty$ for $\mu \in \bcal{O}\union\{\eta^{\Phi}: \Phi\in\bcal{K}_{p}\}$. Suppose also that Assumption \ref{Hyp-Transition-kernel-absolute-continuity} and either conditions $(\mathcal{S}1)$-$(\mathcal{S}4)$ or $(\mathcal{W}1)$-$(\mathcal{W}4)$ are satisfied.}
If there exists $\widetilde{\mu}\in\bcal{O}_{s}$ satisfying $\tilde{\theta}_{i}<\widetilde{\mu}(g_{i})$
for any $i\in \NN_{m}$ then 
\begin{align}
\sup\big\{ \mu(h) :\mu\in\bcal{O} \hbox{ and } & \mu(g_{i})\geq\tilde{\theta}_{i} \text{ for } i\in \NN_{m} \big\} \nonumber \\
& =  \sup\big\{ \mu(h) :\mu\in\bcal{O}_{s} \hbox{ and } \mu(g_{i})\geq\tilde{\theta}_{i} \text{ for } i\in \NN_{m} \big\}.
\end{align}
\end{proposition}
\textbf{Proof:} 
{Let $\bfrak{R}$ be either $\bcal{O}$ or $\{\eta^{\Phi}: \Phi\in\bcal{K}_{p}\}$. Clearly $\beta\mu_{1}+(1-\beta)\mu_{2}\in\bfrak{R}$ for any $\mu_{1}$, $\mu_{2}$ in $\bfrak{R}$ and $\beta\in [0,1]$.}
Let us define $\ds \mathcal{C}= \union_{\mu\in\bfrak{R}} \{\theta\in\mathbb{R}^p: \mu(g_{i})\geq \theta_{i} \text{ for } i\in \NN_{m} \}$.
$\mathcal{C}$ is clearly a non-empty convex subset of $\mathbb{R}^p$. 
Define the function $\mathcal{V}$ on $\mathcal{C}$ by
$$\mathcal{V}(\theta):= \sup\{ \mu(h) :\mu\in\bfrak{R} \hbox{ and } \mu(g_{i})\geq \theta_{i} \text{ for } i\in \NN_{m} \}.$$
By hypothesis, $\mathcal{V}$ takes values in $\RR$ for any $\theta\in\mathcal{C}$.
Observe that $\mathcal{V}$ is a proper concave on~$\mathcal{C}$.
Indeed, consider $\theta_{1}=(\theta_{1,1},\ldots,\theta_{1,m})$ and $\theta_{2}=(\theta_{2,1},\ldots,\theta_{2,m})$ in $\mathcal{C}$ and $\alpha\in [0,1]$.
For any $\epsilon>0$, there exist $\mu_{j,\epsilon}\in\bfrak{R}$ for $j=1,2$ satisfying 
$\mu_{j,\epsilon}(g_{i})\geq \theta_{j,i}$  and $\mu_{j,\epsilon}(h)\geq \mathcal{V}(\theta_{j})-\epsilon/2$ for  $i\in \NN_{m}$.
Clearly, we have $\big(\beta\mu_{1,\epsilon}+(1-\beta)\mu_{2,\epsilon}\big)(g_{i})\geq \beta\theta_{1,i}+(1-\beta)\theta_{2,i}$ for any $i\in \NN_{m}$.
Therefore,
\[\mathcal{V}(\beta\theta_{1}+(1-\beta)\theta_{2}) \geq \big(\beta\mu_{1,\epsilon}+(1-\beta)\mu_{2,\epsilon}\big)(h)
\geq\beta \mathcal{V}(\theta_{1}) +(1-\beta) \mathcal{V}(\theta_{2})-\epsilon,\]
showing that $\mathcal{V}$ is a proper concave function on $\mathcal{C}$.
Now, $\tilde{\theta}$ is in the interior of $\mathcal{C}$, and so $\mathcal{V}$ is continuous at $\tilde{\theta}$ by Proposition 2.17 in \cite{barbu12} and therefore, we can apply Proposition 2.36 in \cite{barbu12} to claim the existence of $\tilde{\lambda}\in\RR^{m}$ such that, for all $\theta\in\mathcal{C}$,
\begin{eqnarray*}
\mathcal{V}(\theta) & \leq & \mathcal{V}(\tilde{\theta}) + \langle \tilde{\lambda},\theta-\tilde{\theta}\rangle.
\end{eqnarray*}
Remark that $\tilde{\lambda} \leq\0_{m}$ since $\mathcal{V}(\theta)\geq \mathcal{V}(\tilde{\theta})$ for all $\theta\leq \tilde{\theta}$.
Now, fix an arbitrary $\mu\in\bfrak{R}$. Then $(\mu(g_{1}),\cdots,\mu(g_{p}))\in\mathcal{C}$ and so,
$$\mathcal{V}(\tilde{\theta})  \geq \mu \big( h-\langle \tilde{\lambda}, g\rangle\big)+\langle \tilde{\lambda}, \tilde{\theta}\rangle.$$
Therefore,
\begin{equation}\label{equation-proof-theorem-LP-1}
\mathcal{V}(\tilde{\theta}) \geq \sup\{ \mu \big( h-\langle \tilde{\lambda}, g\rangle\big) : \mu\in \bfrak{R} \}+\langle \tilde{\lambda}, \tilde{\theta}\rangle.
\end{equation}
For any $\epsilon>0$, there exists $\mu_{\epsilon}\in\bfrak{R}$ with $\mu_{\epsilon}(g_{i})\geq \tilde{\theta_{i}}$ for any $i\in\NN_{m}$ such that $\mu_{\epsilon}(h)\geq \mathcal{V}(\tilde{\theta})-\epsilon$
implying 
\[ \sup\{ \mu \big( h-\langle \tilde{\lambda}, g\rangle\big) : \mu\in \bfrak{R} \}+\langle \tilde{\lambda}, \tilde{\theta}\rangle
\geq \mu_{\epsilon}(h)-\mu_{\epsilon}\big(\langle \tilde{\lambda},g\rangle \big) +\langle \tilde{\lambda}, \tilde{\theta}\rangle
\geq \mu_{\epsilon}(h)\geq \mathcal{V}(\tilde{\theta})-\epsilon\]
since $\tilde{\lambda}\leq\0_{m}$.
Together with (\ref{equation-proof-theorem-LP-1}), this shows 
\begin{equation}\label{equation-proof-theorem-LP}
\sup\{ \mu(h) :\mu\in\bfrak{R} \hbox{ and } \mu(g_{i})\geq\tilde{\theta}_{i} \text{ for } i\in \NN_{m}\} = \sup\{ \mu \big( h-\langle \tilde{\lambda}, g\rangle\big) : \mu\in \bfrak{R} \}
+\langle \tilde{\lambda}, \tilde{\theta}\rangle.
\end{equation}
Now, we have for $\lambda\leq\0_{m}$,
\begin{align*}
\sup \big\{ \mu \big( h-\langle \lambda, g\rangle\big) : \mu\in \bfrak{R} \big\} & +\langle \lambda, \tilde{\theta}\rangle \nonumber \\
& \geq
\sup \big\{ \mu \big( h-\langle \lambda, g\rangle\big) : \mu\in \bfrak{R} \text{ and } \mu(g_{i})\geq \tilde{\theta}_{i} \text{ for } i\in \NN_{m} \big\} +\langle \lambda, \tilde{\theta}\rangle \nonumber \\
& \geq \sup \big\{ \mu(h) :\mu\in\bfrak{R} \hbox{ and }  \mu(g_{i})\geq \tilde{\theta}_{i} \text{ for } i\in \NN_{m} \big\},
\end{align*}
implying
\begin{align*}
\inf \Big\{ \sup\big\{ \mu \big( h-\langle \lambda, g\rangle\big) : \mu\in \bfrak{R} \big\} +\langle \lambda, \tilde{\theta}\rangle  & : \lambda\leq\0_{m}\Big\} \nonumber \\
& \geq \sup\big\{ \mu(h) :\mu\in\bfrak{R} \hbox{ and } \mu(g_{i})\geq \tilde{\theta}_{i} \text{ for } i\in \NN_{m} \big\},
\end{align*}
and so with (\ref{equation-proof-theorem-LP}) we obtain
\begin{align*}
\inf \Big\{ \sup\big\{ \mu \big( h-\langle \lambda, g\rangle\big) : \mu\in \bfrak{R} \big\} +\langle \lambda, \tilde{\theta}\rangle  & : \lambda\leq\0_{m}\Big\} \nonumber \\
& = \sup\big\{ \mu(h) :\mu\in\bfrak{R} \hbox{ and } \mu(g_{i})\geq \tilde{\theta}_{i} \text{ for } i\in \NN_{m} \big\}.
\end{align*}
Therefore, with $\bfrak{R}=\bcal{O}$
\begin{align}
\inf \Big\{ \sup\big\{ \mu \big( h-\langle \lambda, g\rangle\big) : \mu\in \bcal{O} \big\} & +\langle \lambda, \tilde{\theta}\rangle : \lambda\leq\0_{m}\Big\} \nonumber \\
& = \sup\big\{ \mu(h) :\mu\in\bcal{O} \hbox{ and } \mu(g_{i})\geq \tilde{\theta}_{i} \text{ for } i\in \NN_{m} \big\},
\label{Lagrange-General-1}
\end{align}
{and with $\bfrak{R}=\{\eta^{\Phi}: \Phi\in\bcal{K}_{p}\}$
\begin{align}
\inf \Big\{ \sup\big\{ \eta^{\Phi} \big( h-\langle \lambda, g\rangle\big) : \Phi\in\bcal{K}_{p} \big\} & +\langle \lambda, \tilde{\theta}\rangle : \lambda\leq\0_{m}\Big\} \nonumber \\
& = \sup\big\{ \eta^{\Phi}(h) : \Phi\in\bcal{K}_{p} \hbox{ and } \eta^{\Phi}(g_{i})\geq \tilde{\theta}_{i} \text{ for } i\in \NN_{m} \big\}.
\label{Lagrange-Stationary-1}
\end{align}
%By hypothesis, the right hand side of the previous equation is finite and so,
%\begin{align}
%\inf \Big\{ \sup\big\{ \eta^{\Phi} \big( h-\langle \lambda, g\rangle\big) : \Phi\in\bcal{K}_{p} \big\} & +\langle \lambda, \tilde{\theta}\rangle : \lambda\leq\0_{m} \hbox{ and }
% \sup\big\{ \eta^{\Phi} \big( h-\langle \lambda, g\rangle\big) : \Phi\in\bcal{K}_{p} \big\} <+\infty \Big\} \nonumber \\
%&= \inf \Big\{ \sup\big\{ \eta^{\Phi} \big( h-\langle \lambda, g\rangle\big) : \Phi\in\bcal{K}_{p} \big\}  +\langle \lambda, \tilde{\theta}\rangle : \lambda\leq\0_{m}\Big\} \nonumber \\
%& = \sup\big\{ \eta^{\Phi}(h) : \Phi\in\bcal{K}_{p} \hbox{ and } \eta^{\Phi}(g_{i})\geq \tilde{\theta}_{i} \text{ for } i\in \NN_{m} \big\}.
%\label{Lagrange-Stationary-1}
%\end{align}
Now, for $\lambda\leq\0_{m}$ we have $\sup\Big\{ \eta^{\Phi} \Big( \big( h-\langle \lambda, g\rangle\big)^{+}\Big) : \Phi\in\bcal{K}_{p} \Big\} <+\infty$ by hypothesis and we obtain
from Lemma \ref{strategy-stationary-randomized} and Theorem \ref{Linear-prog=stationary-optimality} that
\begin{align}
\sup\big\{ \eta^{\Phi} \big( h-\langle \lambda, g\rangle\big) : \Phi\in\bcal{K}_{p} \big\}  = \sup\{ \mu \big( h-\langle \lambda, g\rangle\big) : \mu\in \bcal{O}_{s} \}
\label{Lagrange-Stationary-2}
\end{align}
and also,
\begin{align}
\sup\big\{ \eta^{\Phi}(h) : \Phi\in\bcal{K}_{p} \hbox{ and } & \eta^{\Phi}(g_{i})\geq \tilde{\theta}_{i} \text{ for } i\in \NN_{m} \big\} \nonumber \\
& = \sup\big\{ \mu(h) :\mu\in\bcal{O}_{s} \hbox{ and } \mu(g_{i})\geq \tilde{\theta}_{i} \text{ for } i\in \NN_{m} \big\}.
\label{Lagrange-Stationary-3}
\end{align}
Therefore, combining equations \eqref{Lagrange-Stationary-1}-\eqref{Lagrange-Stationary-3} we obtain that
\begin{align}
\inf \Big\{ \sup\big\{ \mu \big( h-\langle \lambda, g\rangle\big) : \mu\in \bcal{O}_{s} \big\} & +\langle \lambda, \tilde{\theta}\rangle : \lambda\leq\0_{m}\Big\} \nonumber \\
& = \sup\big\{ \mu(h) :\mu\in\bcal{O}_{s} \hbox{ and } \mu(g_{i})\geq \tilde{\theta}_{i} \text{ for } i\in \NN_{m} \big\},
\label{Lagrange-General-2}
\end{align}
Moreover, Theorem \ref{Schal-Theorem} can be applied to show that
\begin{align}
\sup\{ \mu \big( h-\langle \lambda, g\rangle\big) : \mu\in \bcal{O} \}  = \sup\{ \mu \big( h-\langle \lambda, g\rangle\big) : \mu\in \bcal{O}_{s} \}.
\label{Lagrange-Stationary-4}
\end{align}
Combining equations \eqref{Lagrange-General-1}, \eqref{Lagrange-General-2} and \eqref{Lagrange-Stationary-4}, we obtain the result.}
\hfill$\Box$

\end{document}